\documentclass[10pt,letterpaper]{article}
\usepackage[top=0.85in,left=2.75in,footskip=0.75in]{geometry}

\usepackage{amsmath,amssymb}

\usepackage{changepage}

\usepackage{textcomp,marvosym}

\usepackage{cite}

\usepackage{nameref,hyperref}

\usepackage[right]{lineno}

\usepackage[nopatch=eqnum]{microtype}
\DisableLigatures[f]{encoding = *, family = * }

\usepackage[table]{xcolor}

\usepackage{array}

\newcolumntype{+}{!{\vrule width 2pt}}

\newlength\savedwidth

\usepackage{setspace}

\raggedright
\usepackage[aboveskip=1pt,labelfont=bf,labelsep=period,justification=raggedright,singlelinecheck=off]{caption}

\makeatletter
\renewcommand{\@biblabel}[1]{\quad#1.}
\makeatother

\usepackage{lastpage,fancyhdr,graphicx}
\usepackage{epstopdf}
\fancyheadoffset[L]{2.25in}
\fancyfootoffset[L]{2.25in}
\usepackage{amsthm}
\usepackage{mathtools}
\usepackage{enumitem}
\usepackage{booktabs}
\usepackage{arydshln}
\newcommand{\rowsep}{\noalign{\vskip 1.5pt\global\savedwidth\arrayrulewidth\global\arrayrulewidth 0.2pt}%
  \hdashline[0.5pt/1.5pt]\noalign{\global\arrayrulewidth\savedwidth\vskip 1.5pt}}
\usepackage{xurl}
\usepackage{algorithm}
\usepackage{algpseudocode}
\algnewcommand{\Break}{\textbf{break}}
\algnewcommand{\Continue}{\textbf{continue}}
\newtheorem{theorem}{Theorem}
\newtheorem{lemma}[theorem]{Lemma}
\newtheorem{proposition}[theorem]{Proposition}
\newtheorem{corollary}[theorem]{Corollary}
\newtheorem{conjecture}[theorem]{Conjecture}
\theoremstyle{definition}
\newtheorem{definition}[theorem]{Definition}
\theoremstyle{remark}
\newtheorem{remark}[theorem]{Remark}

\newcommand{\poly}{\operatorname{poly}}
\newcommand{\KL}{\operatorname{KL}}
\newcommand{\Sol}{\mathcal{S}}
\newcommand{\Ber}{\operatorname{Ber}}

\begin{document}
\begin{flushleft}
{\Large
\textbf{A Las Vegas algorithm for satisfiable SAT formulas, with an application to timetabling}
}
\\[1ex]
Fabio Dias\textsuperscript{1*}
\\
\medskip
\textbf{1} Surrey Business School, University of Surrey, Guildford, United Kingdom
\\
\medskip

* f.dias@surrey.ac.uk

\end{flushleft}

\section*{Abstract}
Many scheduling problems need many different feasible timetables, and a view of how much freedom
each event has. We study a Las Vegas algorithm that samples satisfying assignments of a Boolean
formula one variable at a time: it draws each value from an estimate of its conditional
probability, propagates, and restarts on a contradiction. Every output is a checked solution and
comes with the exact probability of producing it. We prove that a run succeeds with probability at least
$2^{-K}$, where $K$ is the expected Kullback--Leibler divergence between the true and estimated
conditional probabilities along the run of a uniform solution. We also bound the distance of the
outputs from uniform from the samples alone; it is a divergence in the opposite direction to $K$, so
a sampler can succeed more often while being farther from uniform. On generated timetabling
instances counted exactly, belief propagation on the timetabling problem brings $K$ to at most
about one bit, and the estimated counts are within 2\% of the exact ones. On a real
university instance of the International Timetabling Competition 2019, rerunning
only the failed runs delivered all 72 requested timetables with each of two estimators. Belief
propagation succeeds in 31\% of runs. Steered by a Markov
chain that carries a nearly feasible timetable, it succeeds in 65\% of runs, with
outputs at least 300 bits from uniform. There the importance weights are too uneven to estimate the
number of timetables, and we certify a lower bound instead. For random 3-SAT, an
estimator whose $K$ grows logarithmically would give a polynomial-time
Las Vegas algorithm beyond the densities at which any algorithm is proved to succeed; we conjecture
that one exists below the satisfiability threshold.

\clearpage
\newgeometry{top=0.85in,left=1in,right=1in,footskip=0.75in}

\section*{Introduction}\label{sec:intro}

Timetabling and other scheduling problems are usually posed as the search for one feasible
solution, or for one optimal solution~\cite{CDS23}. A planner often needs more: several different
feasible timetables to choose from, alternatives that survive a later change of plan, and a view of
which events the constraints pin down and which have room to move. Finding diverse solutions is an
established question in constraint programming~\cite{HHOW05}. Sampling the satisfying assignments
of a Boolean formula uniformly, or nearly so, has dedicated tools: UniGen samples with
near-uniformity guarantees through random parity constraints~\cite{CMV13}, CMSGen and QuickSampler
trade guarantees for speed~\cite{GSCM21,DLBS18}, and Barbarik tests how close a sampler is to
uniform~\cite{CM19}.

\begin{figure}[!htb]
\centering
\includegraphics{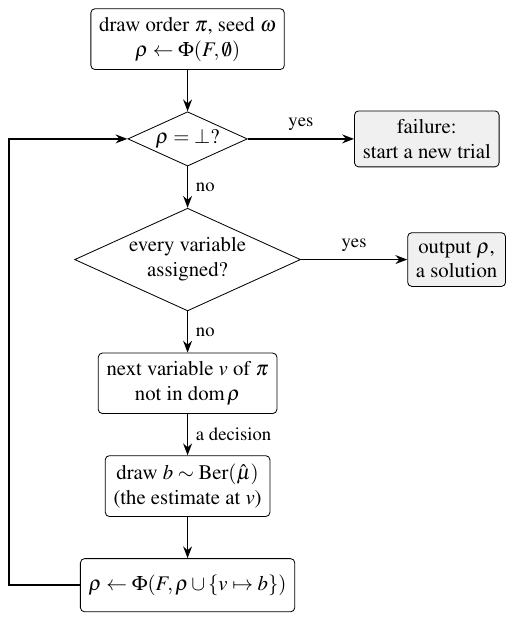}
\caption{\textbf{One trial.}
The closure $\Phi$ (unit propagation in Section~``\nameref{sec:conditional}'') sets
every value it can derive; the estimator $\hat\mu$ (belief propagation) is used only at decisions.
Every output is a solution, so repeated trials form a Las Vegas algorithm.}
\label{fig:trial}
\end{figure}

This paper studies a sampler that also reports exact probabilities. It visits the variables of a
formula in a random order, gives each a value drawn from an estimate of its conditional probability
given the values already set, propagates unit clauses after each step, and restarts on a
contradiction (Fig~\ref{fig:trial}). Every output is a total assignment that falsifies no
clause, so the algorithm never errs, and repeating it is a Las Vegas algorithm. Each output also
comes with the exact probability with which the run produced it: the product of the probabilities
of its random choices. Weighting each output by the inverse of that probability is sequential
importance sampling, which goes back to Knuth's estimator for backtrack trees~\cite{Knu75} and to
Rosenbluth and Rosenbluth~\cite{RR55}. It gives estimates of the number of solutions and
of the fraction of solutions with any property, such as the fraction of feasible timetables that
put a given event in a given slot. For model counting, SampleSearch combines importance sampling
with propagation~\cite{GD11}, ApproxCount and SampleCount sample solutions to estimate
counts~\cite{WS05,GHSS07}, and BPCount uses belief-propagation marginals~\cite{KSS08}.

How often a run succeeds depends on the estimates. Our main theorem (Theorem~\ref{thm:success})
bounds it: a run succeeds with probability at least $2^{-K}$, where $K$ is the expected
Kullback--Leibler divergence between the true and estimated conditional probabilities, summed
along the run of a uniform solution. This is the cost of importance sampling~\cite{CD18} applied
one decision at a time, and it holds for any formula, any sound propagation rule and any estimator.
Exact conditional probabilities give $K = 0$, but each costs a counter call, as in the reduction of
uniform generation to counting~\cite{JVV86}; even the subexponential counter claimed in a recent
preprint~\cite{openai} would give subexponential, not polynomial, time per run. Belief
propagation~\cite{MM09} gives cheap estimates, and $K$ measures what their errors cost.

We apply the algorithm to sampling feasible timetables (Section~``\nameref{sec:timetables}''). The
encoding makes every auxiliary variable a function of the timetable, so each feasible timetable is
exactly one satisfying assignment, and sampling or counting assignments is sampling or counting
timetables. On generated instances with exact counts, belief propagation on the formula turns out
to be confidently wrong several times per run: the auxiliary variables of the encoding form long
loops of deterministic constraints, around which it becomes certain of wrong values. Belief
propagation on the timetabling problem itself, with each event as one variable and each slot's
capacity as one factor, removes these errors. Its $K$ is at most about one bit, its estimated
numbers of timetables are within $2\%$ of the exact counts, and the samples are close to uniform.

We then take the algorithm to a real instance: the university course timetabling problem
\texttt{tg-spr18} of the International Timetabling Competition 2019~\cite{MRM18}, in its reduced
form~\cite{HMSS22}, with $676$ classes, $45\,597$ placements (a time and a room) and about three
million forbidden pairs of placements. Here the number of feasible timetables is far too large to
count, and the experiment answers three questions. Can the Las Vegas algorithm deliver every
timetable it is asked for? Yes: rerunning only the failed runs delivered $72$ of $72$ with each of
two estimators. How close to uniform are its outputs? We bound this from the samples alone
(Proposition~\ref{prop:distance}). And how do two estimators compare? Belief propagation on the
classes succeeds in $31\%$ of runs. Belief
propagation steered by a Markov chain that carries a nearly feasible complete timetable, a
\emph{witness}, succeeds in $65\%$ of runs, but its outputs are at least $300$ bits from uniform. The
two estimators trade success against uniformity, and Section~``\nameref{sec:itc}'' says when to use which.

The same framework gives a conditional result on random 3-SAT, which is secondary here. A random
3-CNF formula with $n$ variables and $\alpha n$ clauses is satisfiable with high probability for
small $\alpha$ and unsatisfiable for large $\alpha$, with a sharp transition~\cite{Fri99} predicted
by the cavity method at $\alpha_s \approx 4.267$~\cite{MPZ02,MMZ06}. The best proved bounds are
$3.52$~\cite{KKL06,HS03} and $4.4898$~\cite{DKMP09}; for large $k$ the threshold of $k$-SAT is
known up to $O(k)$~\cite{AP04} and exists and equals the cavity prediction~\cite{DSS22}. The early
analyses of algorithms are of unit-clause heuristics~\cite{CF86}. The largest density at which an
algorithm is proved to find a solution with probability $1 - o(1)$ is about $3.003$~\cite{FS96};
myopic rules reach $3.26$~\cite{AS00}, and algorithms that use the degrees of literals reach $3.52$
with positive probability~\cite{KKL06,HS03}. For large $k$, the algorithm Fix succeeds up to
$2^k \ln k / k$~\cite{Coj10}, which matches the density at which the solutions shatter into
clusters~\cite{ACO08}. Survey propagation finds solutions close to the threshold in
experiments~\cite{MPZ02,BMZ05}, and with backtracking very close to it~\cite{MPR16}.

 As $\alpha$ grows,
the solutions break into clusters at the dynamical transition $\alpha_d$, and later a few clusters
come to dominate at the condensation transition $\alpha_c$~\cite{KMRSZ07,MRS08}; the connection with
reconstruction on trees is due to M\'ezard and Montanari~\cite{MM06}. For $k = 3$ the two
transitions coincide at $3.86$, set by a local (Kesten--Stigum) instability of the
replica-symmetric solution~\cite{KMRSZ07,BMW00}; the one-step solution is itself unstable below
$\alpha_m \approx 4.15$~\cite{MPR04}. Belief-propagation-guided decimation has been studied
numerically and through its ideal version~\cite{MRS07,RTS09}. For large $k$, it provably fails at
densities $\rho\,2^k/k$~\cite{Coj11,Coj17}, the ideal decimation process has been
analysed~\cite{CP12}, and survey-propagation-guided decimation fails well below the
threshold~\cite{Het16}. Barrier results rule out classes of algorithms at large $k$ or in other
models: sequential local algorithms through the overlap gap property~\cite{GS17}, low-degree
polynomials~\cite{BH21}, stable samplers in shattered phases~\cite{EAG24}, and autoregressive
samplers where the decimated measure has a phase transition~\cite{GDKZ24}. For $k = 3$ we know of no
barrier below $3.86$. Rigorous samplers and approximate counters for random $k$-CNF exist for large
$k$, at densities that have grown from small compared with $2^k$~\cite{GGGY20} to about
$2^{k/3}$~\cite{HWY23} and to $2^k/\mathrm{poly}(k)$~\cite{CLWYY24}; at low density belief
propagation gives the free-entropy density~\cite{MS07}, and Moitra's method counts approximately
for worst-case $k$-CNF of bounded degree~\cite{Moi19}. For $k = 3$ at the densities studied here
none of them applies. Decimation along a sample is pinning, which reduces correlations on
average~\cite{Mon08,RT12,BRS11,CKPZ18}; sampling the Sherrington--Kirkpatrick measure by
algorithmic stochastic localization~\cite{EMS22} uses the same idea.

Our contributions are theoretical (T1--T5) and practical (P1--P3).
\begin{enumerate}[label=(T\arabic*)]
  \item \emph{A success bound for Las Vegas sampling} (Theorem~\ref{thm:success}). Each run
        succeeds with probability at least $2^{-K(F)}$; when the solution is unique, the bound is
        an equality for each fixed order and seed. The bound is the classical cost of sequential
        importance sampling~\cite{Knu75,RR55,CD18}, written for a sound closure and an adaptive,
        possibly randomised, estimator, and it carries over to variables with finite domains
        (Remark~\ref{rem:domains}). Isolating a unique solution~\cite{VV86} squares the accuracy the
        estimates need (Proposition~\ref{prop:isolation}), and frozen variables that propagation
        does not set cost linearly in the error rather than quadratically
        (Proposition~\ref{prop:frozen}).
  \item \emph{A measure of uniformity from the samples alone} (Proposition~\ref{prop:distance}).
        A lower bound on the number of solutions minus the mean of $-\log_2 P$ over the outputs is a
        lower bound on the divergence of the outputs from uniform. It is the divergence in the
        opposite direction to $K$.
  \item \emph{What small $K$ means} (Section~``\nameref{sec:counting}''). An estimator with
        $K = O(\log n)$ yields approximate counting within any fixed factor, and an approximate
        counter of high enough precision yields such an estimator
        (Proposition~\ref{prop:counting}). Such an estimator on every satisfiable 3-CNF exists if
        and only if $\mathrm{NP} = \mathrm{RP}$ (Proposition~\ref{prop:worst}).
  \item \emph{A conditional theorem for random 3-SAT, and a conjecture}
        (Section~``\nameref{sec:conditional}''). If some polynomial-time estimator has $K/\log_2 n$
        bounded in probability over satisfiable formulas at density $\alpha$, a run succeeds with
        polynomially small probability on all but any fixed fraction of formulas. We conjecture
        that such an estimator exists at every density below the satisfiability threshold. For
        belief propagation below $3.86$, we reduce the hypothesis to three lemmas for which no
        proof is known, using a form of the known pinning lemma (Section~``\nameref{sec:towards}'').
  \item \emph{The decimated phase diagram for $k = 3$} (Section~``\nameref{sec:cavity}''), by the method of~\cite{RTS09}, with a positive control at $k = 4$.
\end{enumerate}
\begin{enumerate}[label=(P\arabic*)]
  \item \emph{An encoding and an estimator for timetabling} (Section~``\nameref{sec:timetables}''). The
        encoding has exactly one satisfying assignment per feasible timetable, so the importance
        weights estimate the number of feasible timetables and the share of them that put each
        event in each slot. On generated instances with exact counts, belief propagation on the
        formula is confidently wrong several times per run, around the long deterministic loops
        formed by the encoding's auxiliary variables. Belief propagation on the timetabling problem
        itself removes these errors: $K$ falls to at most about one bit, the estimated numbers of
        timetables are within $2\%$ of the exact counts, and the samples are close to uniform.
  \item \emph{Every requested timetable on a real instance} (Section~``\nameref{sec:itc}''). On a
        competition instance with $676$ classes, the algorithm delivered $72$ of $72$ requested
        timetables with each of two estimators, by rerunning only the failed runs, and every output
        passed a validator. The design fixes the work per step by the structure of the
        problem (arc consistency, singleton probes on the class being decided, and one
        belief-propagation update over a fixed radius), never by a time limit or a convergence test.
  \item \emph{Two estimators, and when to use which} (Section~``\nameref{sec:itc}''). Belief propagation succeeds in $31\%$ of runs. Belief propagation
        steered by a Markov chain's witness succeeds in $65\%$ of runs but is at least $300$ bits
        from uniform. We report the limitation that both leave the importance
        weights too uneven to count on this instance, and give a certified lower bound on the number
        of feasible timetables instead.
\end{enumerate}

We prove no unconditional algorithmic result on random 3-SAT above $3.003$, and no proof of the
conjecture is known for any density. The algorithm is correct on every satisfiable formula; its
polynomial running time is claimed only conditionally, and only for random formulas. The
timetabling experiments use generated instances small enough to count exactly and one real
instance that is not; on the real instance we do not measure $K$, and we compare with no other
sampler. The cavity computation is numerical.

Section~``\nameref{sec:lv-sampling}'' proves the success bound, bounds the distance from uniform, and
states the algorithm with belief propagation. Section~``\nameref{sec:timetables}'' applies it to generated
timetables and Section~``\nameref{sec:itc}'' to a real instance. Section~``\nameref{sec:counting}''
relates $K$ to counting, Section~``\nameref{sec:conditional}'' states the conditional theorem and the
conjecture, Section~``\nameref{sec:towards}'' reduces the hypothesis to three lemmas, and
Section~``\nameref{sec:cavity}'' computes the phase diagram. Section~``\nameref{sec:future}'' concludes.

\section*{Las Vegas sampling with estimated marginals}\label{sec:lv-sampling}

Throughout, $F$ is a CNF formula on variables $V$, $|V| = n$, with a non-empty solution set
$\Sol = \Sol(F) \subseteq \{0,1\}^V$, and $U$ is the uniform distribution on $\Sol$. A
\emph{partial assignment} $\rho$ is a map from a subset $\operatorname{dom}\rho \subseteq V$ to
$\{0,1\}$; $\sigma \supseteq \rho$ means that $\sigma$ agrees with $\rho$ on $\operatorname{dom}\rho$.

\begin{definition}[Sound closure]\label{def:closure}
A \emph{closure} is a polynomial-time map $\Phi$ that sends $(F,\rho)$ either to the symbol
$\bot$ or to a partial assignment $\Phi(F,\rho) \supseteq \rho$, such that
\begin{enumerate}[label=(\roman*)]
  \item if $\Phi(F,\rho) \neq \bot$, every $\sigma \in \Sol$ with $\sigma \supseteq \rho$ satisfies
        $\sigma \supseteq \Phi(F,\rho)$;
  \item $\Phi(F,\rho) = \bot$ only if no $\sigma \in \Sol$ extends $\rho$;
  \item if $\Phi(F,\rho) \neq \bot$ then $\Phi(F,\rho)$ falsifies no clause of $F$.
\end{enumerate}
Unit propagation is a sound closure. So is the following closure, for any constant $D$ at least
the width of the clauses of $F$: derive from $F$ restricted by $\rho$ every polynomial-calculus
consequence of degree at most $D$, set every variable whose value is derived, repeat until nothing
changes, and return $\bot$ if the constant $1$ is derived. It runs in time
$n^{O(D)}$~\cite{CEI96}.
\end{definition}

\begin{definition}[Estimator]\label{def:estimator}
An \emph{estimator} $\hat\mu$ is a polynomial-time rule that, at each decision of the trial below,
returns a number in $[0, 1]$, read as an estimate of $U(\sigma_v = 1 \mid \sigma \supseteq
\rho)$ for the variable $v$ being decided and the current partial assignment $\rho$. It may
depend on $F$, on the order $\pi$ of the trial, on the \emph{history} of the run (the decisions
made so far, with their values, and the partial assignments the closure produced), and on a
random seed $\omega$ drawn once at the start of the trial, independently of everything else. No
clipping is required. If an estimate gives probability $0$ to a value that some remaining solution
takes, the cost $K(F)$ defined below is infinite and Theorem~\ref{thm:success} holds trivially.
Clipping every estimate to $[\eta, 1-\eta]$, with $\eta = \eta(n) > 0$, keeps $K(F)$ finite; exact
marginals keep it finite without clipping. Belief propagation on the residual formula,
warm-started from the messages of the previous step and clipped, is the estimator we have in
mind. An estimator may also use what is known about where the formula comes from, such as the
timetabling instance that it encodes (Section~``\nameref{sec:timetables}'').
\end{definition}

Given $F$, $\Phi$ and $\hat\mu$, one trial is Algorithm~\ref{alg:trial}. A
variable $v$ that reaches line~\ref{alg:trial:decide}, i.e.\ that the closure has not already set,
is a \emph{decision}. The Las Vegas algorithm repeats the trial with fresh coins until it returns a
solution.

\begin{algorithm}[!htbp]
\caption{One trial of Las Vegas sampling with estimated marginals}\label{alg:trial}
\begin{algorithmic}[1]
\Require a CNF formula $F$ on variables $V$, $|V| = n$; a sound closure $\Phi$
  (Definition~\ref{def:closure}); an estimator $\hat\mu$ (Definition~\ref{def:estimator})
\Ensure a solution of $F$, or \textsc{fail}
\State draw a uniformly random order $\pi$ of $V$ and a seed $\omega$
\State $\rho \gets \Phi(F, \emptyset)$
\For{$i = 1, \dots, n$}
  \If{$\rho = \bot$}
    \State \Return \textsc{fail}
  \EndIf
  \State $v \gets \pi(i)$
  \If{$v \notin \operatorname{dom}\rho$} \label{alg:trial:decide}
    \State $p \gets \hat\mu(F, \pi, \omega, \text{history})$
      \Comment{estimate of $U(\sigma_v = 1 \mid \sigma \supseteq \rho)$}
    \State draw $b \sim \Ber(p)$
    \State $\rho \gets \Phi(F, \rho \cup \{v \mapsto b\})$
  \EndIf
\EndFor
\If{$\rho = \bot$}
  \State \Return \textsc{fail}
\EndIf
\State \Return $\rho$ \Comment{a total assignment that falsifies no clause}
\end{algorithmic}
\end{algorithm}

A trial makes at most $n + 1$ calls to $\Phi$ and at most $n$ calls to $\hat\mu$, so it runs in
polynomial time. By Definition~\ref{def:closure}(iii), an output is a total assignment that
falsifies no clause, so \emph{every output is a solution}. Repeating the trial is therefore a Las
Vegas algorithm: it never errs, and only its running time is random.

Write $r = (\pi, \omega)$. Fix $r$ and a solution
$\sigma \in \Sol$, and consider the run of the trial in which every decision $v$ draws
$b = \sigma_v$. By Definition~\ref{def:closure}(i)--(ii), the closure never returns $\bot$ along
this run and only ever assigns values that agree with $\sigma$. So the run outputs $\sigma$. Let
$\mathcal{D}_r(\sigma)$ be its set of decisions, and for $v \in \mathcal{D}_r(\sigma)$ let $\rho_v$ be the partial
assignment just before $v$ is decided. Write
\[
  \mu_v \coloneqq \Ber\bigl(U(\sigma'_v = 1 \mid \sigma' \supseteq \rho_v)\bigr),
  \qquad
  \hat\mu_v \coloneqq \Ber(\text{the estimate at $v$ along this run}),
\]
for the true conditional marginal of $v$ and its estimate, and define
\[
  K_r(F) \coloneqq \mathbb{E}_{\sigma \sim U} \sum_{v \in \mathcal{D}_r(\sigma)} \KL(\mu_v \,\|\, \hat\mu_v),
  \qquad
  K(F) \coloneqq \mathbb{E}_r\, K_r(F),
\]
with $\KL$ the Kullback--Leibler divergence in bits.

\begin{theorem}[Success bound]\label{thm:success}
For every satisfiable $F$, every sound closure $\Phi$ and every estimator $\hat\mu$, one trial
outputs a solution with probability
\[
  p(F) \;\ge\; 2^{-K(F)}.
\]
If $|\Sol(F)| = 1$ then, for each $r$, the success probability $p_r(F)$ of the trial with that
order and seed equals $2^{-K_r(F)}$.
\end{theorem}

\begin{proof}
Fix $r$. For $\sigma \in \Sol$, the trial outputs $\sigma$ exactly when every decision draws
$\sigma$'s value. If it does, the run is the one described above and outputs $\sigma$. If some
decision $v$ draws $1-\sigma_v$, every later assignment extends $\rho_v \cup \{v \mapsto
1-\sigma_v\}$, so the output, if any, differs from $\sigma$. The decisions are made sequentially,
and along the run that follows $\sigma$ each estimate is the one that run computes, so
\begin{equation}\label{eq:P}
  P_r(\sigma) \coloneqq \Pr[\text{output} = \sigma \mid r]
  = \prod_{v \in \mathcal{D}_r(\sigma)} \hat\mu_v(\sigma_v).
\end{equation}
For the uniform measure, Definition~\ref{def:closure}(i) says the events
$\{\sigma' \supseteq \rho\}$ and $\{\sigma' \supseteq \Phi(F,\rho)\}$ coincide under $U$. So the
chain rule along the same run has a factor $1$ at each closure step and
\begin{equation}\label{eq:U}
  U(\sigma) = \prod_{v \in \mathcal{D}_r(\sigma)} \mu_v(\sigma_v).
\end{equation}
Since $U(\sigma) > 0$ on $\Sol$,
\[
  p_r(F) = \sum_{\sigma \in \Sol} P_r(\sigma)
  = \mathbb{E}_{\sigma \sim U}\!\left[\frac{P_r(\sigma)}{U(\sigma)}\right]
  \;\ge\; 2^{\,\mathbb{E}_{\sigma\sim U} \log_2 (P_r(\sigma)/U(\sigma))}
\]
by Jensen's inequality for the concave $\log_2$. By \eqref{eq:P} and \eqref{eq:U},
\[
  \mathbb{E}_{\sigma\sim U} \log_2 \frac{U(\sigma)}{P_r(\sigma)}
  = \mathbb{E}_{\sigma\sim U} \sum_{v \in \mathcal{D}_r(\sigma)} \log_2 \frac{\mu_v(\sigma_v)}{\hat\mu_v(\sigma_v)}
  = K_r(F).
\]
For the last equality, expose $\sigma$ one decision at a time along $\pi$. Whether the $i$-th
variable of $\pi$ is a decision, the assignment $\rho_v$ before it, the history of the run and
hence the estimate are functions of $r$ and of the values of $\sigma$ already exposed. Under $U$,
the event that a solution takes those exposed values coincides with the event that it extends
$\rho_v$, by Definition~\ref{def:closure}(i). So, given the exposed values, $\sigma_v$ has law
$\mu_v$. The conditional expectation of the $i$-th summand is therefore
$\KL(\mu_v\|\hat\mu_v)$, and the tower property gives the claim. Hence
$p_r(F) \ge 2^{-K_r(F)}$. Finally,
$p(F) = \mathbb{E}_r\, p_r(F) \ge \mathbb{E}_r 2^{-K_r(F)} \ge 2^{-K(F)}$ by Jensen's
inequality for the convex $x \mapsto 2^{-x}$.

If $\Sol = \{\xi\}$, the expectation over $U$ is the single point $\xi$, so Jensen's inequality
holds with equality and $p_r(F) = P_r(\xi) = 2^{-K_r(F)}$.
\end{proof}

\begin{corollary}[Las Vegas time]\label{cor:lv}
Independent repetition of the trial outputs a solution after an expected number of at most
$2^{K(F)}$ trials, and never outputs a non-solution. If $K(F) \le c\log_2 n$, the expected time
is polynomial, $n^{c} \cdot \poly(n)$. If $\Pr_F[K(F) > g(n)] \to 0$ over a distribution of
formulas, then almost every formula is solved in expected time $2^{g(n)} \poly(n)$.
\end{corollary}

\begin{proof}
The number of trials until the first success is geometric with mean $1/p(F) \le 2^{K(F)}$.
\end{proof}

\begin{remark}[The bound is only sufficient]\label{rem:sufficient}
A large $K(F)$ does not imply a small $p(F)$. If $F$ has no clauses, every run succeeds, so
$p(F) = 1$. Yet with the estimator $\hat\mu \equiv \eta$, each of the $n$ decisions contributes
$\KL(\Ber(\tfrac12)\|\Ber(\eta)) \to \infty$ as $\eta \to 0$. The bound compares the trial with a
\emph{uniform} sampler, and the trial may succeed by sampling a non-uniform subset of the
solutions. A small $K(F)$ is therefore sufficient for a polynomial success bound, and the theorem
does not claim it is necessary. Proposition~\ref{prop:distance} below measures the other direction,
and Section~``\nameref{sec:itc}'' meets the gap on a real instance and gives a finer bound that accounts for it.
\end{remark}

\begin{remark}[Finite domains, and closures that see the order]\label{rem:domains}
The proof uses only the soundness conditions (i)--(iii) of Definition~\ref{def:closure} and the
product form~\eqref{eq:P}. It therefore holds unchanged in two generalisations used in
Section~``\nameref{sec:itc}''. First, a variable may take values in a finite domain, with $\mu_v$ and
$\hat\mu_v$ laws on that domain and $\KL$ between them; deciding a class of a timetabling problem
among its placements is one such decision, and is the same as a block of binary decisions on the
encoding's primary variables with the same law and the same product of probabilities. Second, the
closure may depend on the order and seed $r$ as well as on $\rho$, for instance by testing the
values of the variable about to be decided, as long as it is sound for each $r$. A randomised
estimator, such as one driven by a Markov chain, is covered by the seed $\omega$: the chain's random
choices form part of it.
\end{remark}

The quantity $K$ needs exact conditional probabilities, which are out of reach on large instances.
What the samples do give is the probability $P_r(\sigma)$ of each output. Write
$b(\sigma) = -\log_2 P_r(\sigma)$ for the \emph{bits} of an output, and let
$P'_r = P_r / p_r(F)$ be the law of the output given that the trial succeeds.

\begin{proposition}[Distance from uniform, from the samples]\label{prop:distance}
Let $Z = |\Sol(F)|$ and fix $r$ with $p_r(F) > 0$.
\begin{enumerate}[label=(\alph*)]
  \item $\log_2 Z - \mathbb{E}_{\sigma\sim P'_r}\, b(\sigma)
        = \KL(P'_r \,\|\, U) - \log_2\bigl(1/p_r(F)\bigr)$. So for any $L \le \log_2 Z$,
        \[
          L - \mathbb{E}_{\sigma\sim P'_r}\, b(\sigma) \;\le\; \KL(P'_r \,\|\, U).
        \]
  \item Over $N$ independent trials, $\hat Z = N^{-1}\sum_i \boldsymbol{1}\{\text{trial $i$
        succeeds}\}\, 2^{b(\sigma_i)}$ has mean $|\operatorname{supp} P_r| \le Z$, with equality
        when $P_r$ is positive on every solution. Hence, for every $\delta \in (0,1)$, with
        probability at least $1 - \delta$,
        \[
          \log_2 Z \;\ge\; \log_2 \hat Z - \log_2(1/\delta).
        \]
\end{enumerate}
\end{proposition}

\begin{proof}
(a) $P'_r(\sigma) = 2^{-b(\sigma)}/p_r(F)$ and $U(\sigma) = 1/Z$ on $\Sol$, so
$\KL(P'_r\|U) = \mathbb{E}_{P'_r} \log_2\bigl(Z\, 2^{-b(\sigma)}/p_r(F)\bigr)
= \log_2 Z - \mathbb{E}_{P'_r} b(\sigma) + \log_2(1/p_r(F))$. The inequality follows from
$p_r(F) \le 1$. (b) Each term has mean $\sum_{\sigma:\,P_r(\sigma)>0} P_r(\sigma)/P_r(\sigma)$,
which is $|\operatorname{supp} P_r|$. By Markov's inequality,
$\Pr[\hat Z \ge Z/\delta] \le \delta$.
\end{proof}

The left side of (a) is computable from a lower bound $L$ on $\log_2 Z$ and the mean bits of the
outputs. When the trial draws a fresh $r$ for each run, the mean over successful runs averages (a)
over $r$, weighted by $p_r(F)$, and still bounds the average divergence of a run's output law from
uniform from below. Both $K$ and $\KL(P'_r\|U)$ are $0$ exactly when every estimate is the true
conditional probability, but they weigh errors differently. $K = \KL(U\|P_r)$ charges the trial for
every solution it rarely produces, so a small $K$ forces a high success rate. $\KL(P'_r\|U)$ charges it
for concentrating on few solutions, and allows a high success rate: when $Z > 1$, a deterministic
solver that always returns the same solution has $p = 1$, the largest possible $\KL(P'\|U)$, namely $\log_2 Z$,
and $K = \infty$. A sampler can therefore raise its success rate by moving towards that extreme,
and Section~``\nameref{sec:itc}'' shows one that does.

Part~(b) needs the support condition. An estimate of exactly $0$ on a value that some remaining
solution takes makes that solution unreachable: $\hat Z$ then estimates the number of reachable
solutions, a lower bound on $Z$. Mixing every estimate with a uniform share $\epsilon$ over the
values that the closure allows, $\hat\mu' = (1-\epsilon)\hat\mu + \epsilon\cdot\text{uniform}$, makes
every solution reachable, at a cost of at most $\log_2\frac{1}{1-\epsilon}$ bits per decision in
$K$. Finally, independent trials amplify exactly: if a trial succeeds with probability $p$, then
rerunning only the failed trials leaves a request unmet after $k$ trials with probability $(1-p)^k$,
which needs the trials to be independent, including any internal state of the estimator.

\begin{proposition}[Isolation squares the accuracy requirement]\label{prop:isolation}
Let $\delta_v \coloneqq |\mu_v(1) - \hat\mu_v(1)|$ be the error at decision $v$.
\begin{enumerate}[label=(\alph*)]
  \item If $|\Sol(F)| = 1$, then $p_r(F) = \prod_{v \in \mathcal{D}_r(\xi)} (1 - \delta_v)$. So, for a
        given $r$, $p_r(F) \ge n^{-c}$ forces $\sum_v \delta_v \le c \ln n$ along that run. With $\Theta(n)$ decisions,
        the average error must be $O(\log n / n)$.
  \item If every decision has a true marginal in $[\gamma, 1-\gamma]$ for a constant
        $\gamma > 0$, and $\delta_v \le \gamma/2$, then
        $\KL(\mu_v\|\hat\mu_v) \le 4\delta_v^2 / (\gamma \ln 2)$. So an average squared error of
        $O(\log n / n)$, i.e.\ errors of order $\sqrt{\log n / n}$, already gives
        $K_r(F) = O(\log n)$.
\end{enumerate}
\end{proposition}

\begin{proof}
(a) With $\Sol = \{\xi\}$ each $\mu_v$ is the point mass at $\xi_v$, so
$\hat\mu_v(\xi_v) = 1 - \delta_v$, and Theorem~\ref{thm:success} gives the product. Then
$-\ln(1-\delta) \ge \delta$ gives $\sum_v \delta_v \le -\ln p_r(F) \le c\ln n$.
(b) In nats, $\KL(P\|Q) \le \ln(1 + \chi^2(P\|Q)) \le \chi^2(P\|Q)$, and for Bernoulli laws
$\chi^2 = (\mu - \hat\mu)^2 / (\hat\mu(1-\hat\mu))$. Here $\gamma \le 1/2$ and
$\hat\mu_v(1) \in [\gamma/2, 1-\gamma/2]$, so $\hat\mu(1-\hat\mu) \ge (\gamma/2)(3/4) \ge \gamma/4$.
Divide by $\ln 2$ for bits.
\end{proof}

\begin{remark}
Valiant--Vazirani isolation~\cite{VV86} leaves a unique solution with probability
$\Omega(1/n)$. By Proposition~\ref{prop:isolation}(a), it makes every error count linearly, where
part~(b) has it count quadratically. It also adds random parity constraints on about $n/2$
variables each, so the isolated formula is no longer a sparse CNF, and belief propagation loses
the local tree structure it relies on. The trial therefore runs on the original formula, without
isolation.
\end{remark}

\begin{proposition}[Frozen variables cost linearly]\label{prop:frozen}
If a decision $v$ is \emph{frozen} in the residual formula ($\mu_v$ is a point mass at $b$)
but not forced by $\Phi$, its contribution is
$\KL(\mu_v\|\hat\mu_v) = -\log_2(1-\delta_v) \ge \delta_v/\ln 2$, linear in the error. A decision
with true marginal in $[\gamma, 1-\gamma]$ and $\delta_v \le \gamma/2$ contributes
$O(\delta_v^2/\gamma)$, quadratic in the error (Proposition~\ref{prop:isolation}(b)).
\end{proposition}

\begin{proof}
For a point mass at $b$, $\KL = \log_2 (1/\hat\mu_v(b)) = -\log_2(1-\delta_v)$, and
$-\ln(1-\delta) \ge \delta$.
\end{proof}

\begin{remark}
Every frozen variable that $\Phi$ derives is set exactly and is not a decision, so it costs
nothing. The linear cost of Proposition~\ref{prop:frozen} falls only on the frozen variables the
closure cannot see. This is the role of the closure in the trial. (A stronger closure also
changes the trajectory and the residuals the estimator sees, so the theorem does not compare
the values of $K$ for two closures.)
\end{remark}

The estimator used in the rest of the paper is belief propagation on the residual formula, with unit propagation as the closure.
Let $\rho$ be the current partial
assignment. A clause is \emph{active} if $\rho$ does not satisfy it. For an active clause $a$ and a
variable $v$ of $a$, write $s_{a,v} \in \{0,1\}$ for the value of $v$ that makes its literal in $a$
true. The message $u_{v\to a} \in [0,1]$ is the probability, in the formula without $a$, that the
literal of $v$ in $a$ is false; for $v \in \operatorname{dom}\rho$ it is fixed at $1$, since an
active clause has only false literals among the assigned variables. Define
\[
  \zeta_{a\to v} \coloneqq \prod_{w \in a,\, w \ne v} u_{w\to a},
  \qquad
  W^{(-a)}_v(x) \coloneqq \prod_{\substack{b \ni v \text{ active},\ b \ne a\\ x \ne s_{b,v}}}
     \bigl(1 - \zeta_{b\to v}\bigr) \quad (x \in \{0,1\}),
\]
and let $W_v(x)$ be the same product over all active $b \ni v$. Here $\zeta_{a\to v}$ is the
probability that the other literals of $a$ are all false. A value $x$ of $v$ that falsifies $v$'s literal in $b$ needs another literal of $b$ to be true,
which has probability $1 - \zeta_{b\to v}$; this is the factor in $W$.

\begin{algorithm}[!htbp]
\caption{Belief propagation on the residual formula}\label{alg:bp}
\begin{algorithmic}[1]
\Require $F$, the partial assignment $\rho$, messages $u$ (warm start), tolerance $\tau$, sweep
  limit $R$
\Ensure updated messages $u$
\For{$s = 1, \dots, R$}
  \State $\Delta \gets 0$
  \For{each $v \notin \operatorname{dom}\rho$, in a fixed order}
    \For{each active clause $a \ni v$}
      \State $u' \gets W^{(-a)}_v(1 - s_{a,v}) \,/\, \bigl(W^{(-a)}_v(0) + W^{(-a)}_v(1)\bigr)$
        \Comment{$1/2$ if both weights are $0$}
      \State $\Delta \gets \max(\Delta, |u' - u_{v\to a}|)$; \quad $u_{v\to a} \gets u'$
        \Comment{in place}
    \EndFor
  \EndFor
  \If{$\Delta < \tau$}
    \State \Return $u$
  \EndIf
\EndFor
\State \Return $u$
\end{algorithmic}
\end{algorithm}

The estimate for a free variable $v$ is
$\hat\mu_v(1) = W_v(1)/(W_v(0) + W_v(1))$, clipped to $[\eta, 1 - \eta]$.

\begin{algorithm}[!htbp]
\caption{The Las Vegas algorithm with belief propagation}\label{alg:lvbp}
\begin{algorithmic}[1]
\Require a CNF formula $F$ on $n$ variables; tolerance $\tau$; sweep limits $R_0$, $R$;
  clipping $\eta$
\Ensure a solution of $F$ (the loop runs until one is found)
\State $u^0 \gets$ Algorithm~\ref{alg:bp} on $(F, \emptyset)$ from $u \equiv 1/2$, with limit $R_0$
  \Comment{shared by all trials}
\Loop \Comment{one iteration is one trial}
  \State $\rho \gets \mathrm{UP}(F, \emptyset)$; \quad $u \gets u^0$
  \While{$\rho \ne \bot$ and some clause of $F$ is active}
    \State $v \gets$ a uniformly random variable not in $\operatorname{dom}\rho$
    \State $p \gets \min\bigl(\max\bigl(W_v(1)/(W_v(0) + W_v(1)),\ \eta\bigr),\ 1 - \eta\bigr)$
    \State draw $b \sim \Ber(p)$
    \State $\rho \gets \mathrm{UP}(F, \rho \cup \{v \mapsto b\})$
      \Comment{unit propagation; $\bot$ on a conflict}
    \If{$\rho \ne \bot$}
      \State $u \gets$ Algorithm~\ref{alg:bp} on $(F, \rho)$ from $u$, with limit $R$
    \EndIf
  \EndWhile
  \If{$\rho \ne \bot$}
    \State set every variable not in $\operatorname{dom}\rho$ to $0$
      \Comment{each appears in no active clause}
    \State \Return $\rho$
  \EndIf
\EndLoop
\end{algorithmic}
\end{algorithm}

Picking a uniformly random unassigned variable at each step has the same law as following a
uniformly random order (Algorithm~\ref{alg:trial}). A trial of Algorithm~\ref{alg:lvbp} stops once
every clause is satisfied. In Algorithm~\ref{alg:trial} the remaining variables would be decisions
with $\mu_v = \hat\mu_v = 1/2$, which succeed and cost nothing in $K$, so the two have the same
success probability and the same $K$.

\section*{Sampling feasible timetables}\label{sec:timetables}

\subsection*{The problem and its encoding}

A timetabling instance has $E$ events, $T$ time slots and $R$ rooms in each slot, with rooms
interchangeable. Every event takes exactly one slot; two events that share a student never take the
same slot; each event may take only its allowed slots (a teacher's or a room's availability); and
at most $R$ events take any slot. A feasible timetable assigns a slot to every event and meets
these hard constraints. Given an instance, a planner may want many different feasible timetables;
the \emph{flexibility map}, the share of feasible timetables that put event $e$ in slot $t$ for
every $e$ and $t$, which shows the events the constraints pin down and those with room to move; and
the number of feasible timetables.

The encoding has a primary variable $y_{e,t}$ for ``event $e$ takes slot $t$''. Exactly one slot
per event is written as a chain: auxiliary variables $r_{e,t}$ with
$r_{e,1} \leftrightarrow y_{e,1}$, $r_{e,t} \leftrightarrow r_{e,t-1} \vee y_{e,t}$ and
$y_{e,t} \rightarrow \neg r_{e,t-1}$ for $t > 1$, and the unit clause $r_{e,T}$. A conflict between
events $a$ and $b$ gives the clauses $\neg y_{a,t} \vee \neg y_{b,t}$, and an unavailable slot the
unit clause $\neg y_{e,t}$. At most $R$ events per slot is a sequential counter~\cite{Sin05} whose
auxiliary variables are defined by equivalences: $s_{t,i,j}$, ``at least $j$ of the first $i$ events
take slot $t$'', satisfies $s_{t,i,j} \leftrightarrow s_{t,i-1,j} \vee (y_{i,t} \wedge
s_{t,i-1,j-1})$, and the clause $\neg s_{t,E,R+1}$ closes it. Every auxiliary variable is
therefore a function of the primary ones, so each feasible timetable is exactly one satisfying
assignment. The uniform law over satisfying assignments is the uniform law over feasible
timetables, and the number of satisfying assignments is the number of feasible timetables. Once
every primary variable is set, unit propagation sets every auxiliary variable, so a run decides only
primary variables: its order is a random order of the $y_{e,t}$, and Theorem~\ref{thm:success}
applies as stated. On small instances we checked, by enumerating every assignment of one slot per
event, that the formula's models are exactly the feasible timetables.

From $N$ runs, of which the successful ones output timetables $\sigma_i$ with probabilities
$P(\sigma_i)$, the number of feasible timetables is estimated by
$\hat Z = N^{-1}\sum_i 1/P(\sigma_i)$. The flexibility map is estimated by the weighted share
$\sum_i [\sigma_i(e) = t]/P(\sigma_i) \big/ \sum_i 1/P(\sigma_i)$.

\subsection*{Two estimators}

The first estimator is belief propagation on the formula (Algorithm~\ref{alg:bp}), applied to the
encoding as to any formula. The second runs belief propagation on the timetabling problem itself.
Each event is one variable, whose domain is the set of slots it may still take given the values set
so far. Each conflict is a factor ``different slots''. Each slot's capacity is one factor ``at most
$R$ events''; its message to event $e$ is the ratio of the probabilities that at most $R - 1$ and
at most $R$ of the other events take the slot, computed exactly from the other events' probabilities
of taking it by a Poisson-binomial dynamic program in $O(ER)$ steps. The estimate of $y_{e,t}$ is
event $e$'s belief in slot $t$. Messages are damped by a factor $0.5$ and iterated to a largest
change below $10^{-4}$, for at most $1000$ sweeps at the start and $100$ after each step, warm
started. The estimator depends on the instance and on the history of the run, and runs in
polynomial time, so it is admissible (Definition~\ref{def:estimator}). The closure is unit
propagation on the encoding, and every output is checked against the formula.

\subsection*{The experiment}

A seeded generator draws the instances. Each of $10$ students takes $3$ of $E = 10$ events chosen
uniformly, and two events conflict when a student takes both; each event may take each slot with
probability $0.8$, keeping at least one slot; there are $R = 3$ rooms. The number of slots
$T \in \{5, 6, 7\}$ is the slack: the same seed gives the same events, conflicts and availability
for every $T$, so a larger $T$ relaxes the same instance. Three seeds give nine instances, with
between $28$ and $2\,677\,680$ feasible timetables, counted exactly by backtracking. One larger
instance has $E = 30$ events, $30$ students, $T = 10$ slots and $R = 4$ rooms; it is too large to
count.

On the six instances with at most $10^5$ feasible timetables we measured $K$ with exact
conditional probabilities: a path follows a run but draws each decision from the exact conditional
probability, obtained by counting the timetables that extend the current values, and adds the
divergence from the estimate. Each $K$ is the mean of $20$ paths ($400$ on the tightest instance).
A decision is \emph{gross} when the estimate, before any clipping, is below $10^{-3}$ or above
$1 - 10^{-3}$ while the exact conditional probability is between $0.05$ and $0.95$. Each sampling
run made up to $20\,000$ runs of the algorithm within $140$ seconds, on one shared Linux server. The
programs, the run scripts and the raw outputs of every experiment in this paper are
available from the author on request.

\subsection*{Results}

Table~\ref{tab:tt-est} compares the two estimators on the six instances where $K$ was measured.
Belief propagation on the formula needs its estimates clipped to $[0.1, 0.9]$. Without clipping,
its estimated numbers of timetables are $0.001$ to $0.14$ of the exact ones; with a pairwise
encoding of one slot per event (an at-least-one clause and a clause for every pair of slots), it
was even so confident on some instances that every run output the same timetable, which the chain
encoding removes. With clipping it is still grossly wrong $3$ to $13$ times per path. At those decisions, the clauses that push its
estimate to near certainty involve the chain of the event's slot variables in $88\%$ to $97\%$ of
cases and the capacity counters in $77\%$ to $87\%$, and conflict clauses in only $19\%$ to $31\%$;
in about $80\%$ of them the certainty comes from belief-propagation messages rather than from
values already set. Each primary variable lies on one event chain and one slot counter, so the
auxiliary variables form a grid of long loops of deterministic constraints, on which belief
propagation on the formula becomes certain of wrong values. We checked the implementation against
exact marginals on $3000$ random formulas whose factor graphs are trees, where belief propagation
is exact: the largest difference was $2\times10^{-16}$. Belief propagation on the timetabling
problem has none of these loops. Its gross decisions fall to at most $0.05$ per path, $K$ falls from
$3.6$--$8.4$ bits to $0.09$--$1.0$ bits, and its runs need no clipping.

\begin{table}[!ht]
\centering
\small
\caption{\textbf{The two estimators on the instances with at most $10^5$ feasible timetables.}}
\label{tab:tt-est}
\begin{tabular}{lrcccccccc}
\toprule
 & & \multicolumn{4}{c}{belief propagation on the formula} & \multicolumn{4}{c}{on the timetabling problem} \\
\cmidrule(lr){3-6}\cmidrule(lr){7-10}
instance & $Z$ & $K$ & gross & $\hat p$ & $\hat Z/Z$ & $K$ & gross & $\hat p$ & $\hat Z/Z$ \\
\midrule
$T = 5$, seed 1 & $1740$ & $6.4$ & $7.8$ & $0.80$ & $0.69$ & $0.41$ & $0$ & $0.89$ & $0.984$ \\
\rowsep
$T = 5$, seed 2 & $28$ & $3.6$ & $3.0$ & $0.41$ & $0.98$ & $1.00$ & $0$ & $0.56$ & $1.004$ \\
\rowsep
$T = 5$, seed 3 & $21\,888$ & $7.4$ & $10.8$ & $0.95$ & $0.81$ & $0.17$ & $0.05$ & $0.98$ & $0.982$ \\
\rowsep
$T = 6$, seed 1 & $33\,523$ & $7.8$ & $10.9$ & $0.95$ & $0.44$ & $0.09$ & $0$ & $0.97$ & $0.994$ \\
\rowsep
$T = 6$, seed 2 & $2076$ & $7.1$ & $8.4$ & $0.84$ & $0.90$ & $0.36$ & $0.05$ & $0.92$ & $0.994$ \\
\rowsep
$T = 7$, seed 2 & $62\,228$ & $8.4$ & $12.8$ & $0.97$ & $0.66$ & $0.11$ & $0.05$ & $0.98$ & $0.993$ \\
\bottomrule
\end{tabular}
\begin{flushleft}$Z$: number
of feasible timetables. $K$: cost in bits, mean over exact paths. Gross: gross decisions per path.
$\hat p$: success rate per run. $\hat Z/Z$: estimated over exact number of timetables. Belief
propagation on the formula is clipped to $[0.1, 0.9]$. On the timetabling problem, $K$ is measured
with the estimates clipped to $[0.01, 0.99]$, and $\hat p$ and $\hat Z/Z$ come from unclipped runs.
\end{flushleft}
\end{table}

Table~\ref{tab:tt-all} gives the results of belief propagation on the timetabling problem on all
nine instances. The estimated numbers of timetables are within $2\%$ of the exact counts, and the
weighted flexibility map is within $0.2$ to $0.3$ percentage points of the exact one on average.
The raw shares of the samples, without weights, are already close (within $0.2$ to $3.3$ points),
so the samples are close to uniform: the mean of $-\log_2 P(\sigma)$ over the outputs matches
$\log_2 Z$ to within $0.1$ bits on every instance except the tightest. Theorem~\ref{thm:success}
holds on every instance; on the tightest, with $28$ timetables, it is nearly tight, with
$2^{-K} = 0.50$ against a success rate of $0.56$ over about $72\,000$ runs. With belief propagation
on the formula the bound held as well, but loosely: success rates of $0.41$ to $0.97$ against
$2^{-K}$ of $0.003$ to $0.09$. On the instance with $30$ events, $554$ of $555$ runs succeeded, at
about four runs per second, and the mean of $-\log_2 P(\sigma)$ was $71.5$ bits.

\begin{table}[!ht]
\centering
\small
\caption{\textbf{Belief propagation on the timetabling problem, without clipping, on all nine instances.}}
\label{tab:tt-all}
\begin{tabular}{lrcccccc}
\toprule
instance & $Z$ & $\log_2 Z$ & $\hat p$ & mean $-\log_2 P$ & $\hat Z/Z$ & map error, raw & weighted \\
\midrule
$T = 5$, seed 1 & $1740$ & $10.77$ & $0.892$ & $10.75$ & $0.984$ & $0.010$ & $0.002$ \\
\rowsep
$T = 5$, seed 2 & $28$ & $4.81$ & $0.564$ & $5.44$ & $1.004$ & $0.033$ & $0.003$ \\
\rowsep
$T = 5$, seed 3 & $21\,888$ & $14.42$ & $0.978$ & $14.34$ & $0.982$ & $0.003$ & $0.002$ \\
\rowsep
$T = 6$, seed 1 & $33\,523$ & $15.03$ & $0.974$ & $14.98$ & $0.994$ & $0.006$ & $0.003$ \\
\rowsep
$T = 6$, seed 2 & $2076$ & $11.02$ & $0.916$ & $10.99$ & $0.994$ & $0.014$ & $0.002$ \\
\rowsep
$T = 6$, seed 3 & $417\,312$ & $18.67$ & $0.995$ & $18.64$ & $0.997$ & $0.003$ & $0.002$ \\
\rowsep
$T = 7$, seed 1 & $272\,424$ & $18.06$ & $0.993$ & $18.02$ & $1.000$ & $0.004$ & $0.002$ \\
\rowsep
$T = 7$, seed 2 & $62\,228$ & $15.93$ & $0.985$ & $15.87$ & $0.993$ & $0.007$ & $0.002$ \\
\rowsep
$T = 7$, seed 3 & $2\,677\,680$ & $21.35$ & $0.999$ & $21.34$ & $1.002$ & $0.002$ & $0.002$ \\
\bottomrule
\end{tabular}
\begin{flushleft}Mean $-\log_2 P$: over the successful runs. Map error: mean absolute error of the flexibility
map, from the raw shares of the samples and from the importance-weighted shares.
\end{flushleft}
\end{table}

These instances are generated and small enough to count exactly. Section~``\nameref{sec:itc}'' takes the
algorithm to a real instance that is not.

\section*{A real instance: the International Timetabling Competition 2019}\label{sec:itc}

\subsection*{The instance and the model}

The International Timetabling Competition 2019~\cite{MRM18} poses university course timetabling
problems taken from real universities. Each class must get a time, made of weeks, days, a start and
a length, and a room, from lists given in the instance; rooms have unavailable times; no two classes
may use a room at overlapping times; and the instance lists distribution constraints between
classes, some required and some penalised. We use the reduced form of the competition
instances~\cite{HMSS22}, in the same format, and the instance \texttt{tg-spr18}. It has $676$
classes, $24$ rooms and $45\,767$ placements (a time and a room), of which $45\,597$ remain once
placements in unavailable room-times are removed. Its required distribution constraints are all
between pairs of classes: $67$ SameAttendees, $23$ SameDays, $17$ NotOverlap, $16$ DifferentDays,
$12$ SameTime, $4$ Precedence and $1$ SameRoom, with the competition's semantics; SameAttendees
uses the travel times between rooms. We keep the hard constraints only. The penalties and the
assignment of students are optimisation criteria, and a feasible timetable is a placement for every
class that meets every hard constraint.

Every hard constraint then forbids pairs of placements of two classes. There are $2\,969\,369$
forbidden pairs: $782\,526$ come from room clashes and $2\,329\,950$ from SameAttendees, whose $67$
constraints have up to $174$ classes each; a pair can come from more than one constraint. The
encoding is that of Section~``\nameref{sec:timetables}'': a primary variable for each placement of each
class, exactly one placement per class as a chain, and one binary clause per forbidden pair. It has
$91\,194$ variables and $3\,151\,081$ clauses, and each feasible timetable is exactly one
satisfying assignment. A validator re-reads the competition file and checks every
output against the hard constraints from their definitions, not from the encoder's pair lists.

The number of feasible timetables is far too large to count. The experiment therefore asks whether
the algorithm delivers every timetable requested of it, how close to uniform its outputs are
(Proposition~\ref{prop:distance}), and how two estimators compare.

\subsection*{The trial}

Decisions are made on classes (Remark~\ref{rem:domains}). A step picks an undecided class $c$
uniformly at random, which gives the same law as a uniformly random order. It then runs a
\emph{singleton probe} on each placement $o$ still allowed for $c$: arc consistency on the current
domains, with $c$ restricted to $o$. Every placement whose probe fails is in no feasible timetable
that agrees with the current domains, and is removed. The step then draws a placement of $c$ from
the estimate, renormalised over the placements left, and applies the closure: arc consistency, which
removes every placement that has no compatible placement left in some other class, repeated until
nothing changes, together with unit propagation on the encoding. Arc consistency only removes
placements that no feasible completion uses, and the probes depend only on the state and on $c$,
so both are sound closures in the sense of Remark~\ref{rem:domains}. A run fails only when the
probes leave a class with no placement. At the start, the closure removes $4\,278$ of the $45\,597$
placements.

\emph{Belief propagation on the classes.} Each class is a variable over its allowed placements, and
each pair of classes with forbidden pairs between them is a factor; there are $13\,905$ factors. The
message from class $a$ to class $b$ at placement $o'$ is one minus the probability, in $a$'s
cavity belief, of the placements of $a$ that are forbidden with $o'$. Before the first run, belief
propagation runs to a largest change below $1/N$, for at most $N$ sweeps, where $N = 676$ is the
number of classes; on this instance it stopped by its tolerance after $100$ sweeps, in $10$ seconds.
After each decision, every factor within distance $2$ of the classes whose domains changed is
updated once, warm-started. On this instance that neighbourhood is nearly the whole graph, so a
decision costs about one sweep.

\emph{Belief propagation steered by a witness.} A Markov chain runs on complete timetables within
the current domains, feasible or not. A timetable with $v$ violated forbidden pairs has weight
$\lambda^{v}$, with $\lambda = 1/N$, and a step of the chain picks a class uniformly and redraws its
placement from its conditional given all the other classes (a heat-bath step). Since
$\lambda > 0$, the chain needs no feasible timetable to start from and cannot be trapped by a hard
constraint. At a decision on class $c$, the chain runs one sweep of $N$ steps, and its estimate for
$c$ is $c$'s conditional, given the other classes, averaged over the $N$ states of the sweep. The
steered estimate of a placement is proportional to its belief times the chain's estimate. After a
decision and its closure, every class whose current placement has left its domain is redrawn. The
chain's current timetable is the \emph{witness}: when it violates no pair, it is a feasible
completion of the current partial timetable. Each run starts its own chain from a uniformly random
timetable within the domains and runs it for $N$ sweeps first, so that runs are independent.

\emph{Support.} A belief can be exactly $0$ on a placement that some feasible timetable uses. Both
estimators therefore mix in a uniform share $\epsilon = 1/N$ over the allowed placements, so that
every feasible timetable can be produced and $\hat Z$ is unbiased (Proposition~\ref{prop:distance}).

Every parameter is set by the instance: the tolerance, $\lambda$ and $\epsilon$ are $1/N$, the
chain runs one sweep per decision and $N$ sweeps at the start, and the radius of the update is fixed.
The work per step is set by the structure of the problem, not by a time limit, a count that the
problem does not give, or a test of whether belief propagation has converged. 

\subsection*{How the design was reached}

With unit propagation as the closure and belief propagation run
to convergence after each decision, a run took about a minute, and later several: after a decision,
loopy belief propagation needed $10^5$ to $3.5\times10^6$ factor updates to settle again. Fixing the
work per step as above and adding arc consistency made a run take about $10$ seconds, but only
$7.6\%$ of runs succeeded ($168$ of $2\,208$). The failures were spread over the whole run, at a
roughly constant rate of about $0.4\%$ per decision, at states from which no feasible timetable
remained but which arc consistency did not detect. The singleton probes raised the success rate to
$34.6\%$ ($64$ of $185$), and every remaining failure was a class left without a placement by its
probes. The witness raised it to $67.5\%$. The last step added the support mixing and the
independent chains, and asked each process for a fixed number of timetables.

One idea did not work. Sequential Monte Carlo over the trial, with populations of $32$ partial runs
resampled after each decision and weighted towards the belief-propagation (Bethe) estimate of the
number of completions, finished no timetable in four populations, while $32$ independent trials
would have finished about $13$. The weight of a decision is the ratio of the Bethe counts after and
before it, divided by the probability of the draw. It rewards a draw that belief propagation rated
unlikely unless the Bethe count falls by the same factor, and the two disagree most near the dead
ends that arc consistency does not see. Resampling steered the population into them.

\subsection*{Results}

Each estimator ran as $36$ independent processes, side by side on one server with $80$ cores. Each
process asked for two timetables and reran its failed runs until it had them. The method has no
time limit; only the experiment as a whole had one, which was not reached.
Table~\ref{tab:itc} gives the results.

\begin{table}[!ht]
\centering
\small
\caption{\textbf{The two estimators on \texttt{tg-spr18}, with arc consistency, the singleton probes, the
support mixing and independent runs.}}
\label{tab:itc}
\begin{tabular}{>{\raggedright\arraybackslash}p{94mm}cc}
\toprule
 & belief propagation & steered \\
\midrule
timetables delivered (processes complete) & $72$ of $72$ ($36$ of $36$) & $72$ of $72$ ($36$ of $36$) \\
\rowsep
outputs that pass the validator & $72$ & $72$ \\
\rowsep
runs & $233$ & $111$ \\
\rowsep
success rate $\hat p$ ($95\%$ confidence interval) & $0.309$ ($0.25$, $0.37$) & $0.649$ ($0.55$, $0.73$) \\
\rowsep
runs per timetable (fewest, most per process) & $3.24$ ($2$, $18$) & $1.54$ ($2$, $8$) \\
\rowsep
decisions with an allowed placement at probability $0$, before and after mixing & $180$, $0$ & $99$, $0$ \\
\rowsep
$\log_2 \hat Z$ & $2173.2$ & $1880.9$ \\
\rowsep
mean bits of the outputs (interquartile range) & $2160.7$ ($16.2$) & $1838.5$ ($26.4$) \\
\rowsep
distance from uniform, at least, against the certified $L = 2165.6$ & $4.9$ & $327.1$ \\
\rowsep
distance from uniform, at least, against the largest estimate, $2177.8$ & $17.1$ & $339.3$ \\
\rowsep
effective sample size & $2.2$ of $72$ & $1.0$ of $72$ \\
\rowsep
classes that differ between two outputs, same and different processes & $0.808$, $0.806$ & $0.801$, $0.801$ \\
\rowsep
time for the arm & $20$ min & $16$ min \\
\bottomrule
\end{tabular}
\begin{flushleft}Bits: $-\log_2 P$ of an output. The distance from uniform is
the lower bound of Proposition~\ref{prop:distance}(a), against the certified lower bound $L$ on
$\log_2 Z$ and against the largest $\log_2 \hat Z$ of all our runs on this instance, which is not
certified. Effective sample size: $(\sum_i w_i)^2/\sum_i w_i^2$ for the importance weights
$w_i = 2^{b(\sigma_i)}$. Success rate: $\hat p = 72/T$ over the $T$ runs of the arm; the confidence
interval is the exact one for a sequence of independent runs stopped at its $72$nd success.
\end{flushleft}
\end{table}

\emph{Every timetable was delivered.} Both arms delivered $72$ of $72$ timetables, and every output
passed the validator. The number of runs per timetable equals $1/\hat p$ by
construction; it is within $13\%$ of the $1/p$ predicted from the earlier success rates for belief
propagation, and within $5\%$ for the steered estimator. The unluckiest process needed $18$ runs
with belief propagation and $8$ with the steered estimator. Because runs are independent, rerunning
only the failed runs leaves a request unmet after $k$ runs with probability $(1-p)^k$; at the
measured rates, $19$ and $7$ runs per timetable give $99.9\%$. A run took about three minutes of one
core with $72$ processes sharing the machine.

\emph{The support gap was real.} Before the mixing, $180$ decisions of belief propagation and $99$ of
the steered estimator gave probability $0$ to a placement that the closure allowed. Without the
mixing, some feasible timetables could not have been produced, and $\hat Z$ would have estimated
only the number of reachable ones. After the mixing there were none.

\emph{A certified lower bound on the number of feasible timetables.} Each process stops at its
second success, and the runs of $36$ independent processes, taken in turn, form one sequence of
independent runs stopped at its $72$nd success. With $S$ successes in $T$ runs, the estimate
$\hat Z = \frac{S-1}{T-1}\cdot\frac{1}{S}\sum_{i} 2^{b(\sigma_i)}$, summed over the successes, is
unbiased, because $(S-1)/(T-1)$ is unbiased for the success probability under this stopping
rule~\cite{Hal45} and is independent of the bits of the successes. By
Proposition~\ref{prop:distance}(b), applied to both arms with $\delta = 0.005$ each, with
probability at least $0.99$,
\[
  \log_2 Z \;\ge\; 2173.2 - \log_2 200 \;=\; 2165.6 .
\]

\emph{Distance from uniform.} Against this certified bound, the outputs of belief propagation are at
least $4.9$ bits from uniform and those of the steered estimator at least $327$ bits
(Proposition~\ref{prop:distance}(a)). Against the largest estimate of $\log_2 Z$ from all our runs,
$2177.8$, which is not certified, the figures are $17.1$ and $339.3$ bits.

\emph{Diversity.} Two outputs differ on about $80\%$ of classes with either estimator, whether they
come from the same process or from different ones. So the outputs are diverse, and runs in one
process are as different from each other as runs in different processes.

\emph{The limitation: uneven weights.} The bits of the outputs spread over an interquartile range of
$16$ and $26$ bits, so their importance weights span many orders of magnitude, and the effective
sample size is $2.2$ and $1.0$ of $72$. An importance-weighted estimate on this instance, of the
number of timetables or of the share of timetables that put a class in a placement, rests on one or
two outputs. The flexibility map of Section~``\nameref{sec:timetables}'' cannot be estimated here, and the
count is reported only as a certified lower bound. On the generated instances, where $K$ was at
most about one bit, the weights were even; here they are not, and we did not measure $K$.

\subsection*{Why the steered estimator succeeds more often, and when to use which}

The steered estimator doubles the success rate and moves its outputs more than $300$ bits farther
from uniform. Theorem~\ref{thm:success} is not contradicted: the bound $p \ge 2^{-K}$ is sufficient,
not necessary (Remark~\ref{rem:sufficient}), and the distance in Table~\ref{tab:itc} is the
divergence in the opposite direction to $K$ (Proposition~\ref{prop:distance}). A high success rate
with a large divergence from uniform is the behaviour of a sampler moving towards a solver: it
concentrates on a region where it rarely meets a dead end.

The proof of Theorem~\ref{thm:success} gives a finer bound that accounts for this. For any
probability distribution $\nu$ on $\Sol$ and any $r$, by Jensen's inequality,
\[
  p_r(F) \;\ge\; \sum_{\sigma \in \operatorname{supp}\nu} P_r(\sigma)
  = \mathbb{E}_{\sigma\sim\nu}\!\left[\frac{P_r(\sigma)}{\nu(\sigma)}\right]
  \;\ge\; 2^{-\KL(\nu \,\|\, P_r)},
\]
and $\nu = U$ gives Theorem~\ref{thm:success}. A trial is then charged only for the solutions it aims
at.  On the generated instances of
Section~``\nameref{sec:timetables}'', belief propagation on the formula succeeded far more often
than $2^{-K}$. On this
instance, the steered estimator succeeds twice as often as belief propagation while its outputs are
hundreds of bits farther from uniform. In both cases a $\nu$ concentrated on the solutions the
trial actually reaches would explain the success rate where $K$ does not.

We traced the witness in an earlier run of the steered estimator, with
$104$ timetables, before the support mixing and with one chain per process carried from run to run.
Along successful runs the witness had a median of $0$ violated pairs at a decision. While it
violates none and lies within the domains, the current partial timetable has a feasible completion,
so it is not a dead end; belief propagation carries no such certificate. In failed runs the median
was $2$, throughout the run and not only near the failure, and over the last $10$ decisions of a
failed run the witness had $2.0$ violated pairs on average, against $1.0$ for successful runs at the
same steps. The draws do not copy the witness: they equal its placement in only $20\%$ of decisions.
But the chain's estimate is a class's conditional given a complete timetable, which in a tightly
constrained problem is narrow. 

When to use which:
\begin{itemize}
  \item \emph{Belief propagation} when the outputs are to stand for the feasible timetables as a
        whole: alternatives meant to be representative, statistics over the timetables, and any
        importance-weighted estimate. Its outputs are the closer to uniform of the two. It needs
        about twice as many runs per timetable.
  \item \emph{The steered estimator} when the aim is many different feasible timetables quickly,
        and how they represent the whole does not matter: it needs about half as many runs per
        timetable, and its outputs still differ on about $80\%$ of classes. Its outputs are far from
        uniform, and it should not be used for counting.
\end{itemize}
Both give every output with its exact probability, a certified lower bound on the number of
feasible timetables, and every requested timetable by rerunning only the failed runs.

\section*{From counting to search}\label{sec:counting}

This section connects the success bound to counting, to the counting algorithm
of~\cite{openai} and to our earlier note~\cite{DL26}.

Let the estimator return the exact conditional marginal,
$\hat\mu_v = \mu_v$. Then $K(F) = 0$, and Theorem~\ref{thm:success} gives $p(F) = 1$: no decision
ever takes a value that no remaining solution has, so the closure never returns $\bot$. Each exact
marginal is a ratio of two counts, of the solutions extending $\rho \cup \{v \mapsto 1\}$ and
$\rho$, so a counting algorithm with running time $T_{\mathrm{count}}(n)$ gives a Las Vegas
algorithm that succeeds on every trial and runs in time $O(n\,T_{\mathrm{count}}(n))$. This is the reduction of Jerrum, Valiant and
Vazirani~\cite{JVV86} from uniform generation to counting. (Search alone needs less: by
self-reducibility it needs only to know whether a count is positive, but still one call per
decision.) With the bound claimed
in~\cite{openai}, it would sample a uniform solution in time $2^{O(n \log\log n/\sqrt{\log n})}$. By
the constants extracted in~\cite{DL26}, that running time beats brute force on no formula of bounded
density and feasible size.

Theorem~\ref{thm:success} places this at one end of a range. At the other end, an estimator that
runs in polynomial time makes each trial polynomial, and pays for its errors with a success
probability of at least $2^{-K}$. The rest of the paper studies that end.

Between the two ends sits approximate counting.
Proposition~\ref{prop:counting} relates it to small $K$ in both directions, at different
precisions. Write $Z = |\Sol(F)|$. For an order and seed $r$, let $P_r$ be the output law of the
trial, so $P_r(\sigma) = \prod_{v \in \mathcal{D}_r(\sigma)} \hat\mu_v(\sigma_v)$ by~\eqref{eq:P}; the trial
can compute it along its run, and it is positive on every solution whenever $K_r(F)$ is finite,
for instance under clipping.

\begin{proposition}[Small $K$ and approximate counting]\label{prop:counting}
Let $F$ be satisfiable.
\begin{enumerate}[label=(\alph*)]
  \item If $K(F) \le C\log_2 n$, with $C$ known, for a trial that runs in polynomial time, then for
        every fixed $\epsilon, \beta \in (0,1)$ a randomized algorithm running in time
        $(e\,n^{C})^{O(1/\epsilon^2)}\poly(n)\log(1/\beta)$ outputs $\hat Z$ with
        $\Pr\bigl[|\hat Z/Z - 1| > \epsilon\bigr] \le \beta$.
  \item Suppose an algorithm approximates, within a factor in $[1 - \epsilon_n, 1 + \epsilon_n]$
        and with probability at least $1 - n^{-3}$, the number of solutions of every formula $F$
        restricted by a partial assignment, and always returns $0$ when that number is $0$. Then
        the estimator $\hat\mu_v(1) = \hat c_1/(\hat c_0 + \hat c_1)$, where $\hat c_b$ is the
        approximate count for $v \mapsto b$ (and $\hat\mu_v(1) = 1/2$ if both counts are $0$,
        which happens only when a call fails), clipped with $\eta = n^{-3}$, has
        $K(F) \le 5 n \epsilon_n + o(1)$ for $n^{-2} \le \epsilon_n \le 1/3$. With
        $\epsilon_n = \log_2 n / n$ this is $K(F) = O(\log n)$, at the cost of $2n$ counts per
        trial.
\end{enumerate}
\end{proposition}

\begin{proof}
(a) By Markov's inequality, $K_r(F) \le 4C\log_2 n$ for at least $3/4$ of the pairs $r$. Fix
such an $r$ and draw $N$ independent trials with it. The estimator
$\hat Z = N^{-1}\sum_i \boldsymbol{1}\{\text{trial $i$ succeeds}\}/P_r(\sigma_i)$ is importance
sampling of the uniform measure $U$ with proposal $P_r$, and $\hat Z/Z$ is the estimate of
$\mathbb{E}_U[1] = 1$. $P_r$ is a sub-probability on $\Sol$; extend it to a probability by an atom
``fail'' that carries the probability of failure, has mass $0$ under $U$ and contributes $0$ to
$\hat Z$. In nats, $L \coloneqq \KL(U\|P_r) = K_r(F)\ln 2$ by the proof of
Theorem~\ref{thm:success}. Chatterjee and Diaconis~\cite[Thm.~1.1]{CD18} give, for $N = e^{L+t}$,
\[
  \mathbb{E}\,|\hat Z/Z - 1| \le e^{-t/4} + 2\sqrt{\Pr_{\sigma \sim U}[Y > L + t/2]},
  \qquad Y = \ln\bigl(U(\sigma)/P_r(\sigma)\bigr).
\]
Since $\ln x \le x$ and $\mathbb{E}_U[P_r/U] \le 1$, the negative part of $Y$ has mean at most
$1$, so $\mathbb{E}_U[Y_+] \le L + 1$ and $\Pr[Y > L + t/2] \le (L+1)/(L + t/2)$. Taking
$t = 800(L+1)/\epsilon^2$ makes the right-hand side at most $\epsilon/4$, so
$\Pr[|\hat Z/Z - 1| > \epsilon] \le 1/4$ by Markov's inequality. The algorithm does not know $L$,
so it takes $N = e^{L' + t'}$ with $L' = 4C \ln n \ge L$ and $t' = 800(L'+1)/\epsilon^2$; the bound
of~\cite{CD18} only improves as $t$ grows, and $N = (e\,n^{C})^{O(1/\epsilon^2)}$. A good $r$ and a
good estimate occur together with probability at least $9/16$; the median of $O(\log(1/\beta))$
independent repetitions, each with a fresh $r$, is within the factor except with probability
$\beta$.

(b) Let $\hat c_b$ approximate $c_b$, the number of solutions extending $\rho \cup \{v \mapsto b\}$.
When both calls succeed, $\hat c_b = 0$ exactly when $c_b = 0$, and otherwise
$\hat c_b / c_b \in [1-\epsilon_n, 1+\epsilon_n]$, so the unclipped estimate satisfies
$|\ln(\hat\mu_v(b)/\mu_v(b))| \le \ln\frac{1+\epsilon_n}{1-\epsilon_n} \le 3\epsilon_n$ for every
$b$ with $\mu_v(b) > 0$. Clipping changes $\hat\mu_v(b)$ by at most a factor $1 - \eta$ on such
$b$. Hence $\KL(\mu_v\|\hat\mu_v) \le \max_b \log_2(\mu_v(b)/\hat\mu_v(b)) \le
(3\epsilon_n + 2\eta)/\ln 2 < 5\epsilon_n$ per decision, for at most $n$ decisions. A failed call
happens with probability at most $2n \cdot n^{-3}$ per trial, and costs at most
$\log_2(1/\eta) = 3\log_2 n$ bits per decision, which adds $o(1)$ to $K(F)$.
\end{proof}

The two parts are not an equivalence. Part~(a) gives a fixed-factor approximation of the count of
$F$ itself. Part~(b) needs precision $\log n/n$ on every restricted formula met along the run,
which (a) does not provide. What they show is that small $K$ and approximate counting are close in
strength: an estimator with $K = O(\log n)$ is at least as strong as approximate counting within
any fixed factor, and a sufficiently precise counter gives such an estimator. For random $k$-CNF,
polynomial-time approximate counting is known for large $k$, up to densities of order
$2^k/\mathrm{poly}(k)$~\cite{GGGY20,HWY23,CLWYY24}. For $k = 3$ it is easy at very low density:
below $\alpha = 1/6$ the formula falls into components of size $O(\log n)$ that can be counted
exactly. At low density, belief propagation also gives the free-entropy density~\cite{MS07}, which
is weaker than a fixed-factor count. Near the densities studied here we know of no
approximate counter. A proof of the hypothesis of Section~``\nameref{sec:conditional}'' at such a density
would therefore also give one.

\section*{A conditional theorem for random 3-SAT}\label{sec:conditional}

The random formula $F_n(\alpha)$ has $m = \lfloor \alpha n \rfloor$ clauses, each drawn
independently and uniformly from the $8\binom{n}{3}$ clauses on three distinct variables.
$F_n(\alpha)$ is satisfiable with high probability for $\alpha < 3.52$~\cite{KKL06,HS03} and
unsatisfiable with high probability for $\alpha > 4.4898$~\cite{DKMP09}. In between, the
satisfiability threshold exists as a sharp, possibly $n$-dependent, sequence~\cite{Fri99}, but
its location is not proved. We work with $F_n^{\mathrm{sat}}(\alpha)$, the law of $F_n(\alpha)$
conditioned on being satisfiable.

Below $3.52$, satisfiability is proved through algorithms that set variables one at a
time~\cite{KKL06,HS03}. They succeed with probability bounded away from $0$, and the sharp
threshold~\cite{Fri99} turns that into satisfiability with high probability. The largest density
at which an algorithm is proved to find a solution with probability $1 - o(1)$ is about $3.003$,
for the generalised unit clause rule with limited backtracking~\cite[Thm.~1.2]{FS96}. So $3.003$
is the benchmark for search, and $3.52$ the one for satisfiability. By \emph{Target~A} we mean a
Las Vegas algorithm that, at some density $\alpha > 3.003$, runs in polynomial time per trial
and, for every $\epsilon > 0$, succeeds in each trial with probability at least $1/p_\epsilon(n)$
on a $1 - \epsilon$ fraction of $F \sim F_n^{\mathrm{sat}}(\alpha)$, for some polynomial
$p_\epsilon$. Target~A is a different notion from success with high probability, which is what
the algorithm of~\cite{FS96} gives below $3.003$, and which implies Target~A there. Whether the
algorithms that reach $3.52$ with positive probability~\cite{KKL06,HS03} also meet Target~A, that
is, whether rerun with fresh coins they succeed with polynomially small probability on almost every
formula, is not established as far as we know.

The closure $\Phi$ is unit propagation. The estimator $\hat\mu$ is belief
propagation on the residual formula. It starts from uniform messages, is warm-started from the
previous step's messages, is iterated until the largest message change is below a tolerance
$\tau$ or $R_{\max}$ sweeps have run, and is clipped to $[\eta, 1-\eta]$ with $\eta = n^{-3}$.
It depends on the history of the run, which Definition~\ref{def:estimator} allows. The algorithm
the conjecture is about takes $\tau = 1/n$ and $R_0 = R = n$ in Algorithm~\ref{alg:lvbp} (so
$R_{\max} = n$). A sweep that recomputes each message from scratch costs
$O(\sum_v \deg(v)^2)$, so each trial runs in time $O(n \cdot R_{\max} \cdot \sum_v \deg(v)^2)$,
which is $O(n^3)$ in expectation on random formulas of bounded density and polynomial on every
formula. A tolerance that stays fixed as $n$ grows would leave a
truncation error of order $\tau$ on each decision, and so a cost of order $n\tau^2$ in $K$, linear
in $n$.
Algorithms~\ref{alg:bp} and~\ref{alg:lvbp} in Section~``\nameref{sec:lv-sampling}'' state the belief
propagation and the whole algorithm.

\begin{definition}[Hypothesis $\mathrm{H}_T(\alpha)$]
Let $T$ be a trial (Algorithm~\ref{alg:trial}) whose closure and estimator run in polynomial time,
and $K_T(F)$ the cost $K(F)$ of Section~``\nameref{sec:lv-sampling}'' for this trial.
$\mathrm{H}_T(\alpha)$ says that $K_T(F)/\log_2 n$ is bounded in probability
over $F \sim F_n^{\mathrm{sat}}(\alpha)$: for every $\epsilon > 0$ there are $C_\epsilon$ and
$n_\epsilon$ with $\Pr[K_T(F) \le C_\epsilon \log_2 n] \ge 1 - \epsilon$ for all
$n \ge n_\epsilon$. We write $\mathrm{H}_{\mathrm{BP}}(\alpha)$ when $T$ is
Algorithm~\ref{alg:lvbp}.
\end{definition}

\begin{theorem}[Target~A, conditional]\label{thm:conditional}
Fix $\alpha$ and a trial $T$ as above.
\begin{enumerate}[label=(\alph*)]
  \item Assume $\mathrm{H}_T(\alpha)$. Repeating $T$ is a Las Vegas algorithm that never outputs a
        non-solution. For every $\epsilon > 0$ and $n \ge n_\epsilon$, on a $1 - \epsilon$ fraction
        of $F \sim F_n^{\mathrm{sat}}(\alpha)$, each trial succeeds with probability at least
        $n^{-C_\epsilon}$, so the expected running time is at most $n^{C_\epsilon}\poly(n)$. So $T$
        achieves Target~A at density $\alpha$ if $\alpha > 3.003$.
  \item Assume only $\mathbb{E}_{F \sim F_n^{\mathrm{sat}}(\alpha)} K_T(F) \le C\log_2 n$. Then
        for every $g(n) \to \infty$, on a $1 - o(1)$ fraction of formulas each trial succeeds with
        probability at least $n^{-C g(n)}$.
\end{enumerate}
\end{theorem}

\begin{proof}
Both follow from Theorem~\ref{thm:success} and Corollary~\ref{cor:lv}.
(a) On the event $K_T(F) \le C_\epsilon\log_2 n$, $p(F) \ge n^{-C_\epsilon}$.
(b) By Markov's inequality, $\Pr[K_T(F) > g(n)\,C\log_2 n] \le 1/g(n)$.
\end{proof}

\begin{remark}[Why not one constant for almost every formula]\label{rem:strong}
A stronger hypothesis would ask for $\Pr[K_T(F) \le C\log_2 n] \to 1$ with a single $C$. For
belief propagation we do not expect it, by the following heuristic. Bounded-size structures on
which belief propagation is grossly wrong, such as short cycles, should have Poisson counts with
$n$-free means, and a decision on one of them can cost of order $\log_2(1/\eta) = \Theta(\log n)$
bits. If such a cost occurs with probability bounded away from $0$ per structure, then
$\Pr[K > C\log_2 n]$ stays bounded away from $0$ for every fixed $C$. The bounded-in-probability form is what
Proposition~\ref{prop:towards} delivers. It gives polynomial time on all but any fixed fraction of
formulas, with a degree that grows as the fraction shrinks.
\end{remark}

Theorem~\ref{thm:conditional} is proved. What is not proved is its
hypothesis. The conjecture is about the existence of a good estimator, not about belief
propagation in particular.

\begin{conjecture}\label{conj:main}
For every $\alpha$ with $\liminf_{n} \Pr[F_n(\alpha) \text{ is satisfiable}] > 0$, that is,
roughly, at every density below the satisfiability threshold (predicted at
$4.267$~\cite{MPZ02,MMZ06}; the threshold may depend on $n$~\cite{Fri99}), there
is a trial $T$, with a sound closure and an estimator that run in polynomial time, such that
$\mathrm{H}_T(\alpha)$ holds.
\end{conjecture}

Our candidate below $3.86$ is Algorithm~\ref{alg:lvbp}: unit propagation and belief propagation,
with $\tau = 1/n$ and $R_0 = R = n$. Above $3.86$ we propose no specific estimator (see point~(b)
below); estimators that account for clusters, such as survey propagation and its backtracking
version~\cite{MPZ02,BMZ05,MPR16}, are natural candidates there. Measurements of one estimator
estimate the behaviour of that estimator only. A growth rate measured for belief propagation would not show that a better estimator does not exist. By Proposition~\ref{prop:counting}, the conjecture
sits between two counting statements. It implies polynomial-time approximate counting within any
fixed factor on all but any fixed fraction of formulas. A sufficiently precise approximate counter
on the restricted formulas would imply it.

\begin{corollary}[What a proof of the conjecture gives]\label{cor:main}
If $\mathrm{H}_T(\alpha)$ is proved for some $\alpha$ and some trial $T$ as above, then repeating
$T$ is a Las Vegas algorithm for random 3-SAT at density $\alpha$. It never outputs a non-solution
on any satisfiable 3-CNF, and for every $\epsilon > 0$, on a $1 - \epsilon$ fraction of satisfiable formulas, it finds a
solution in expected time $n^{C_\epsilon}\poly(n)$. If Conjecture~\ref{conj:main} is proved, this
holds at every density below the satisfiability threshold, far above the largest density, $3.003$,
at which an algorithm is proved to find solutions with high probability~\cite{FS96}. By
Proposition~\ref{prop:counting}(a), it also gives approximate counting within any fixed factor on
the same formulas, in time polynomial for each fixed factor. This is not a fully polynomial
scheme.
\end{corollary}

\begin{proof}
Theorem~\ref{thm:conditional}(a), and Proposition~\ref{prop:counting}(a) for the last part. When
$T$ is Algorithm~\ref{alg:lvbp}, which stops once every clause is satisfied and sets the $f$
remaining variables to $0$, the importance weight of an output is $2^{f}/P$, where $P$ is the
product of the probabilities the run used: those $f$ variables appear in no active clause, so all
$2^{f}$ completions are solutions, and Algorithm~\ref{alg:trial} would set each of them with
probability $1/2$.
\end{proof}

Our evidence is below $3.86$, and the
reduction of Section~``\nameref{sec:towards}'' covers only that range. We state the conjecture up to the
threshold because nothing known rules it out. None of the four points below is positive evidence
above $3.86$; they explain why what we know does not argue against it.
\begin{enumerate}[label=(\alph*)]
  \item \emph{It allows polynomial decay.} The conjecture allows the success probability of a trial
        to fall with $n$, as $n^{-C}$ with a constant that depends on the density. A success rate
        that falls with $n$ is not evidence against it; only a super-polynomial decay would be. A
        measured success rate $p$ gives only the lower bound $K \ge \log_2(1/p)$.
  \item \emph{The obstacle above $3.86$ concerns belief propagation, not the best estimator.} The
        cavity method predicts that above $3.86$ the measure is condensed before any variable is
        decimated (at decimated fraction $\theta = 0$, in the notation of
        Section~``\nameref{sec:cavity}''), so that belief
        propagation's first estimates are off by $\Theta(1)$ on a positive fraction of variables.
        Over the early part of the run this could cost $\Theta(n)$ for belief propagation, which
        is why we name it as a candidate only below $3.86$. An estimator that accounts for the
        clusters, as survey propagation does~\cite{MPZ02,BMZ05}, need not pay this cost in the
        same way; we have not measured one.
  \item \emph{For $k = 3$, the planted computation loses reconstruction at small $\theta$ above
        $3.86$.} In the cavity computation of Section~``\nameref{sec:cavity}'', revealing a fraction
        $\theta = 0.0005$ of the planted solution removes reconstruction at
        $\alpha \in \{3.90, 4.00, 4.20\}$. We repeated this with two more seeds at $3.90$ and
        $4.20$, and with a population of $160\,000$ at $3.90$; every repeat gave the same result. Above condensation the planted
        solution is not a uniform one, so this describes the decimation of a planted solution;
        whether it carries over to uniform solutions, and so to the trial, is open. We also do not
        understand why it happens: at $k = 4$, reconstruction survives such fractions above
        condensation.
  \item \emph{Nothing known rules it out.} Survey propagation finds solutions close to the
        threshold in experiments~\cite{MPZ02,BMZ05}, and we know of no barrier result for $k = 3$
        below the threshold.
\end{enumerate}

Correctness holds on every satisfiable 3-CNF: the algorithm never outputs a non-solution, and on
every satisfiable formula its expected time is at most $2^{K_T(F)}\poly(n)$
(Theorem~\ref{thm:success}, Corollary~\ref{cor:lv}). Only the polynomial bound is about random
formulas. Proposition~\ref{prop:worst} gives the exact status of a polynomial bound on every
formula.

\begin{proposition}[Logarithmic $K$ on every formula is equivalent to $\mathrm{NP} = \mathrm{RP}$]
\label{prop:worst}
\begin{enumerate}[label=(\alph*)]
  \item There is a trial $T$, with a sound closure and an estimator that run in polynomial time,
        and a constant $C$ such that $K_T(F) \le C\log_2 n$ for every satisfiable 3-CNF formula
        $F$ on $n$ variables, if and only if $\mathrm{NP} = \mathrm{RP}$.
  \item If a Las Vegas algorithm runs in polynomial expected time on every satisfiable 3-CNF
        formula, even when it is only required to work under the promise that its input is
        satisfiable, then $\mathrm{NP} = \mathrm{RP}$.
\end{enumerate}
\end{proposition}

\begin{proof}
(a), only if. Run $T$ for at most $2n^C$ trials and answer ``satisfiable'' if some trial outputs an assignment.
If $F$ is satisfiable, each trial succeeds with probability at least $2^{-K_T(F)} \ge n^{-C}$ by
Theorem~\ref{thm:success}, so all $2n^C$ trials fail with probability at most
$(1 - n^{-C})^{2n^C} \le e^{-2} < 1/2$. If $F$ is unsatisfiable, no trial outputs anything, since
every output is a solution. The running time is polynomial, so this is a one-sided-error
algorithm for SAT, and $\mathrm{SAT} \in \mathrm{RP}$ gives $\mathrm{NP} = \mathrm{RP}$.

(a), if. Approximating the number of satisfying assignments of a CNF formula within a factor
$1 \pm \epsilon$ is possible in randomized polynomial time, polynomial in $n$ and $1/\epsilon$,
with an NP oracle~\cite{Sto83,JVV86}. If $\mathrm{NP} = \mathrm{RP}$, then
$\mathrm{NP} \subseteq \mathrm{BPP}$ and the oracle can be replaced by a randomized algorithm, so
the approximation runs in randomized polynomial time; repeating it and taking the median makes its
failure probability at most $n^{-3}$. It applies to every formula restricted by a partial
assignment, since such a formula is again a CNF formula. To make it return $0$ exactly when the
count is $0$, first run the one-sided-error algorithm for SAT that $\mathrm{NP} = \mathrm{RP}$
provides, which never accepts an unsatisfiable formula, and return $0$ if it rejects. Proposition~\ref{prop:counting}(b) with
$\epsilon_n = \log_2 n/n$ then gives an estimator, clipped with $\eta = n^{-3}$, with
$K_T(F) \le 5\log_2 n + o(1)$ on every satisfiable formula, using $2n$ approximate counts per trial.
The estimator uses internal randomness, which Definition~\ref{def:estimator} allows.

(b) For a Las Vegas algorithm with expected time $p(n)$ on satisfiable inputs, run it for $2p(n)$
steps and check any output against $F$. By Markov's inequality it answers within that time with
probability at least $1/2$ on satisfiable inputs, and on unsatisfiable inputs no output passes the
check, whatever the algorithm does there. So the promise does not help.
\end{proof}

Conjecture~\ref{conj:main} is the random-formula form of Proposition~\ref{prop:worst}(a). By Proposition~\ref{prop:counting}(a), an estimator whose $K$ is logarithmic on every satisfiable formula would also give approximate
counting of every 3-CNF within a fixed factor, which decides SAT. The conjecture as stated concerns
random formulas, and the $2^{K_T(F)}\poly(n)$ bound is what holds on every formula.
Section~``\nameref{sec:towards}'' reduces the hypothesis for belief propagation below $3.86$ to three
lemmas.

\begin{remark}[What refuting the conjecture would take]\label{rem:refute}
If $\mathrm{NP} = \mathrm{RP}$, the estimator in the proof of Proposition~\ref{prop:worst}(a) has
$K_T(F) \le 5\log_2 n + o(1)$ on every satisfiable formula, so $\mathrm{H}_T(\alpha)$ holds at every
density and Conjecture~\ref{conj:main} holds. Hence a proof that, at some density below the
satisfiability threshold, no trial with a polynomial-time closure and estimator satisfies
$\mathrm{H}_T(\alpha)$ would prove $\mathrm{NP} \ne \mathrm{RP}$, and with it
$\mathrm{P} \ne \mathrm{NP}$. Showing that one particular estimator, such as belief propagation,
fails would not. Measurements can support or weaken a candidate estimator, but they cannot refute the conjecture.
\end{remark}

\begin{remark}[Why the hypothesis is on $K$ and not on a $\chi^2$ sum]\label{rem:chi2}
Since $\KL \le \ln(1+\chi^2) \le \chi^2$ in nats, the sum $X(F)$ of the
$\chi^2(\mu_v\|\hat\mu_v)$ bounds $K(F) \ln 2$, and a hypothesis on $X$ would also imply
Target~A. The converse fails badly for an estimate that is confident and wrong. For $\mu_v(1) =
1/2$ and $\hat\mu_v(1) = \eta$, the decision costs $\tfrac12\log_2\frac{1}{4\eta(1-\eta)}
\approx \tfrac32 \log_2 n$ bits in $K$ but $\chi^2 \approx 1/(4\eta) = n^3/4$. If such decisions occur a bounded number of times per run, the $\chi^2$ sum of the trial is of
order $1/\eta$, which is $n^3$ for the clipping $\eta = n^{-3}$, while $K$ stays $O(\log n)$. In general, a decision costs at most
$\log_2(1/\eta) = 3\log_2 n$ bits, and a frozen decision on which the estimate is saturated at the
correct value, $1 - \eta$, costs $\log_2 \frac{1}{1-\eta} \le 2\eta$ bits. With up to $n$ such
decisions, these cost at most $2n\eta$ bits in total. So any clipping $\eta = n^{-c}$ with
$c \ge 1$ keeps a bounded number of confident errors, and the saturated frozen decisions, within
the $O(\log n)$ that the hypothesis allows.
\end{remark}

Write $\alpha^*_{\mathrm{rec}}$ for
the density below which the ideal decimation trajectory never meets reconstruction
(Section~``\nameref{sec:cavity}''); we find $\alpha^*_{\mathrm{rec}} \approx 3.86$. Two kinds of
evidence bear on $\mathrm{H}_{\mathrm{BP}}(\alpha)$ for $\alpha < \alpha^*_{\mathrm{rec}}$, and so
on the conjecture there. Both rest on the same replica-symmetric cavity picture.
\begin{enumerate}[label=(\roman*)]
  \item \emph{Cavity method} (Section~``\nameref{sec:cavity}''). For $k = 3$, decimating along a uniform
        solution keeps the residual formula out of reconstruction at every decimated fraction
        $\theta$ of our grid for $\alpha < \alpha^*_{\mathrm{rec}} \approx 3.86$, with the same
        result for several seeds and population sizes. This is our computation, by the method
        of~\cite{RTS09}, with a positive control at $k = 4$ that reproduces~\cite{RTS09}. The cavity method predicts that in this phase belief propagation
        gives asymptotically exact marginals on typical variables~\cite{KMRSZ07}.
  \item \emph{Scaling} (a heuristic). Loop corrections predict an $O(\chi_{\mathrm{SG}}/n)$ error
        on a typical decision, $\chi_{\mathrm{SG}}$ being the spin-glass susceptibility, finite
        below $\alpha^*_{\mathrm{rec}}$. The cost is quadratic in such errors
        (Proposition~\ref{prop:isolation}(b)), so typical decisions cost $O(\chi_{\mathrm{SG}}^2/n)$
        in total. Bounded-size structures on which belief propagation is grossly wrong (short
        cycles, small dense subformulas) have Poisson counts with $n$-free means, and each costs
        at most $\log_2(1/\eta) = O(\log n)$ bits (Remark~\ref{rem:chi2}). Together these give
        $K(F) = O_P(\log n)$, which is $\mathrm{H}_{\mathrm{BP}}(\alpha)$. This is a large-$n$
        statement: at
        $\alpha = 3.80$ the correlation length is about $60$ levels of the tree
        (Section~``\nameref{sec:cavity}''), more than the radius up to which a formula of a few thousand
        variables looks like a tree.
        Section~``\nameref{sec:towards}'' turns this heuristic into a reduction to three lemmas with no
        known proof.
\end{enumerate}

Below $\alpha = 1/6$, $\mathrm{H}_T(\alpha)$ holds trivially: with high probability the formula
splits into components of size $O(\log n)$, and exact counting on each component gives an
estimator with $K = 0$. No proof of $\mathrm{H}_T(\alpha)$ is known for any trial at any density
$\alpha \ge 1/6$. (Below $3.003$, Target~A holds through the algorithm of~\cite{FS96}, not through
a bound on $K$.) For belief
propagation, a proof would need rigorous control of it along an adaptive decimation of random
3-SAT. The rigorous
analyses of decimation we know of~\cite{Coj11,CP12} are for large $k$, and~\cite{Coj11} proves
that belief-propagation-guided decimation fails there below the best algorithmic density.

\section*{Towards a proof of the hypothesis}\label{sec:towards}

This section reduces $\mathrm{H}_{\mathrm{BP}}(\alpha)$, and so the conjecture at $\alpha$, for
$\alpha < 3.86$, to a known entropy bound, a bound on the number of short cycles, and three lemmas
for which no proof is known. We also prove a first-order version of the first lemma; the reduction
does not use it, but it shows which normalisation the lemma needs. By Proposition~\ref{prop:counting}, the conclusion
would also give approximate counting, so the open lemmas carry the difficulty of that problem.

Draw $\sigma \sim U$ and a uniformly random order $\pi$, and
reveal $\sigma$ one variable at a time along $\pi$. This is the run in the definition of $K$: the
decisions are the revealed variables that the closure has not already set, and the closure's
values are functions of the revealed ones. Revealing coordinates of a sample is what the
statistical-physics literature calls \emph{pinning}~\cite{Mon08,CKPZ18}, and its basic property is
that pinning removes correlations, on average, at a summable rate. Lemma~\ref{lem:budget} is a form
of the known pinning, or global-correlation-reduction, lemma~\cite{Mon08,RT12,BRS11,CKPZ18}; we
include its short proof in the form we use.

\begin{lemma}[Entropy budget]\label{lem:budget}
Let $S_t = \{\pi(1), \dots, \pi(t)\}$, let $J$ be uniform on $V$ and independent of everything else,
and let
\[
  \kappa_t \coloneqq \mathbb{E}_{\pi, J}\, I\bigl(\sigma_{\pi(t+1)};\ \sigma_J \bigm| \sigma_{S_t}\bigr)
\]
in bits, for $0 \le t < n$. This holds for every satisfiable $F$, indeed for every probability
distribution of $\sigma$ on $\{0,1\}^V$:
\[
  \sum_{t=0}^{n-1} \kappa_t \;\le\; 1 .
\]
In particular, at most $1/\epsilon$ steps have $\kappa_t \ge \epsilon$.
\end{lemma}

\begin{proof}
Fix $\pi$ and $j \in V$. By the chain rule,
$H(\sigma_j \mid \sigma_{S_t}) - H(\sigma_j \mid \sigma_{S_{t+1}}) =
I(\sigma_j; \sigma_{\pi(t+1)} \mid \sigma_{S_t})$, and the sum over $t$ telescopes to
$H(\sigma_j) - H(\sigma_j \mid \sigma) = H(\sigma_j) \le 1$. Average over $\pi$ and $j$.
\end{proof}

Lemma~\ref{lem:budget} bounds the average over $j$, which is
$n^{-1}\sum_j I(\sigma_v;\sigma_j\mid\text{past})$ at the next variable $v$. It says that the
correlation susceptibility of the decided variable, measured in mutual information, is at most
$1$ summed over the run. The average includes $j = v$ itself, which contributes
$H(\sigma_v \mid \text{past})/n$. This is a statement about the run as a whole, and it holds for
every measure, condensed or not; it is not a form of the finiteness of the susceptibility
$\chi_{\mathrm{SG}}$ used in Section~``\nameref{sec:conditional}'', which can fail. The heuristic of
Section~``\nameref{sec:conditional}'' uses $\chi_{\mathrm{SG}}$; the reduction below uses this budget
instead. The work lies in turning it into a bound on the estimator's errors.

Fix a step at which
$v$ is decided, and let $N(v)$ be
the variables that share an active clause with $v$. The \emph{cavity measure} is the law of
$\sigma_{N(v)}$ in the residual formula with $v$ and its clauses removed, given the past. Let
$\tilde\mu_v$ be the belief-propagation update at $v$ computed from the exact marginals of the
cavity measure. This is an oracle, since those marginals are what belief propagation tries to
compute. Call $v$ a \emph{defect} if it lies within a fixed distance $r_0$ of a cycle of length
at most $\ell_0$ in the factor graph of $F$, for constants $r_0, \ell_0$. Write $\nu$ for the
cavity measure, $\bar\nu$ for the product of its marginals, $\mathrm{TC}(v) = \KL(\nu\|\bar\nu)$
for its total correlation (in nats), $A_x$ for the set of neighbour configurations under which
the value $x$ of $v$ satisfies all active clauses of $v$, and
$Z_v = \nu(A_0) + \nu(A_1) \in (0, 2]$ for the cavity expected number of allowed values of $v$.
Then $\mu_v(x) = \nu(A_x)/Z_v$, and $\tilde\mu_v$ is the same with $\bar\nu$ in place of $\nu$.

\begin{lemma}[One step, first order]\label{lem:onestep}
$\KL(\mu_v\|\tilde\mu_v) \le \bigl(2\,\mathrm{TC}(v) + 3\sqrt{\mathrm{TC}(v)}\bigr)/Z_v$, in nats.
The factor $1/Z_v$ cannot be dropped: there are cavity measures with $\mathrm{TC}(v) \to 0$ and
$\KL(\mu_v\|\tilde\mu_v) \to \ln 2$.
\end{lemma}

\begin{proof}
Let $P(x, y) = \nu(y)\boldsymbol{1}\{y \in A_x\}/Z_v$ and $Q(x, y) = \bar\nu(y)\boldsymbol{1}\{y \in A_x\}
/\bar Z_v$, with $\bar Z_v = \bar\nu(A_0) + \bar\nu(A_1)$. Their first marginals are $\mu_v$ and
$\tilde\mu_v$, so $\KL(\mu_v\|\tilde\mu_v) \le \KL(P\|Q) = \mathbb{E}_P[\psi] + \ln(\bar Z_v/Z_v)$,
with $\psi = \ln(\nu/\bar\nu)$. The second marginal of $P$ is $\nu(y)w(y)/Z_v$ with
$w(y) = \#\{x : y \in A_x\} \in \{0,1,2\}$, so $\mathbb{E}_P[\psi] \le 2\,\mathbb{E}_\nu[\psi_+]/Z_v$.
Now $\mathbb{E}_\nu[\psi_+] = \mathrm{TC}(v) + \mathbb{E}_\nu[(\ln(\bar\nu/\nu))_+] \le \mathrm{TC}(v)
+ \mathrm{TV}(\nu,\bar\nu)$, since $(\ln r)_+ \le (r-1)_+$. Also
$\ln(\bar Z_v/Z_v) \le (\bar Z_v - Z_v)/Z_v \le 2\,\mathrm{TV}(\nu,\bar\nu)/Z_v$, since
$w \le 2$. Pinsker's inequality, $\mathrm{TV} \le \sqrt{\mathrm{TC}/2}$, gives the bound. For the
second claim, take the clauses $(v \vee a \vee b)$, $(\neg v \vee c \vee d)$ and
$(\neg v \vee c' \vee d')$. Under $\nu$, let $a, b, c, d$ be independent, with
$\Pr[a \vee b] = \lambda^2$ and $\Pr[c \vee d] = \lambda$, and let $c' = c$ and $d' = d$. Then
$\mu_v(1) = \lambda/(\lambda + \lambda^2) \to 1$, while $\tilde\mu_v(1) = \lambda^2/(\lambda^2 +
\lambda^2) = 1/2$, and $\mathrm{TC}(v)$, the information between $(c,d)$ and its copy, tends to $0$
as $\lambda \to 0$.
\end{proof}

The sum over the run needs a second-order bound, quadratic rather than linear in the cavity
correlations, which Lemma~\ref{lem:onestep} does not give.
\begin{enumerate}[label=(L\arabic*)]
  \item \emph{One step, second order.} $\KL(\mu_v\|\tilde\mu_v) \le C_1\,\mathrm{TC}(v)/Z_v^2$.
        The normalization is needed, by Lemma~\ref{lem:onestep}; the care needed is with
        near-frozen variables, where the bound must stay linear in $\mathrm{TC}$.
  \item \emph{Neighbours look like random pairs.} If $v = \pi(t+1)$, then
        \[
          \mathbb{E}\bigl[\mathrm{TC}(v)\,Z_v^{-2}\,\boldsymbol{1}\{v \text{ is a decision and not a
          defect}\}\bigr] \le C_2\, \mathbb{E}_F\, \kappa_t,
        \]
        with the expectation on the left over $F$, $\sigma$ and $\pi$. In a random formula the neighbours of $v$
        are at distance of order $\log n$ from each other once $v$ is removed. By the symmetry of
        the ensemble, and by quiet planting~\cite{ACO08,KZ09}, which is expected to hold below
        condensation, they should be correlated like random pairs. Coja-Oghlan and Perkins proved
        a qualitative form of this step for replica-symmetric random factor graph models: pinning
        makes the measure nearly symmetric, and near symmetry makes the marginals an approximate
        fixed point of belief propagation~\cite{CP18,CKPZ18}. Our use needs it with rates, along
        an adaptive decimation, and for the correlations of several variables at once, not only of
        pairs.
  \item \emph{Belief propagation finds the oracle's fixed point.} If $v = \pi(t+1)$, the
        estimate $\hat\mu_v$ of Algorithm~\ref{alg:lvbp}, started from uniform messages and
        warm-started along the run, satisfies
        $\mathbb{E}\bigl[\bigl(\KL(\tilde\mu_v\|\hat\mu_v) + \max_b
        \ln^2(\tilde\mu_v(b)/\hat\mu_v(b))\bigr)\boldsymbol{1}\{v \text{ is a decision and not a
        defect}\}\bigr] \le C_3\, \mathbb{E}_F\, \kappa_t$. On the tree, belief propagation from
        uninformative messages reaching the right fixed point is non-reconstruction. The gap $g$ of
        Section~``\nameref{sec:cavity}'' measures, on the tree, the distance between belief propagation
        without and with the planted boundary, so that section's computation is evidence for
        (L3) along the trajectory below $3.86$. (L3) itself concerns belief propagation on the
        finite, loopy formula, warm-started and truncated, which the tree computation does not
        cover.
\end{enumerate}

\begin{proposition}[Conditional on (L1)--(L3)]\label{prop:towards}
Assume (L1)--(L3) at density $\alpha$, and that the number $N_{\mathrm{def}}(F)$ of defect variables of
$F \sim F_n^{\mathrm{sat}}(\alpha)$ is bounded in probability (as it is when short cycles have
Poisson counts with $n$-free means). Then $\mathrm{H}_{\mathrm{BP}}(\alpha)$ holds: for every
$\epsilon > 0$ there is $C_\epsilon$ with $\Pr[K(F) \le C_\epsilon \log_2 n] \ge 1 - \epsilon$
for all large $n$. By Theorem~\ref{thm:conditional}, each trial then succeeds with probability at
least $n^{-C_\epsilon}$ on a $1 - \epsilon$ fraction of formulas.
\end{proposition}

\begin{proof}[Proof, given the lemmas]
For Bernoulli laws $P$, $Q$, $S$, in nats,
\[
  \KL(P\|S) = \KL(P\|Q) + \mathbb{E}_P\ln(Q/S)
  = \KL(P\|Q) + \KL(Q\|S) + \sum_b (P(b) - Q(b)) \ln(Q(b)/S(b)).
\]
With Pinsker's inequality and $2ab \le a^2 + b^2$,
\[
  \KL(P\|S) \le \tfrac32\, \KL(P\|Q) + \KL(Q\|S) + 2 \max_b \ln^2\bigl(Q(b)/S(b)\bigr).
\]
Apply it with $P = \mu_v$, $Q = \tilde\mu_v$, $S = \hat\mu_v$ at every decision that is not a
defect. By (L1) pointwise and (L2)--(L3) in expectation, its expected cost is at most
$C\,\mathbb{E}_F\kappa_t$ nats, with $C = \tfrac32 C_1C_2 + 2C_3$ and $\kappa_t$ in bits. By
Lemma~\ref{lem:budget}, applied to each $F$, the expected total over these decisions is at most
$C$ nats, that is, $C/\ln 2$ bits. By Markov's inequality it is bounded in probability over $F$. A
defect decision costs at most $\log_2(1/\eta) = 3\log_2 n$ bits, because the clipped estimate is
at least $\eta$. There are at most $N_{\mathrm{def}}(F)$ defect decisions. Hence, in bits,
$K(F) \le O_P(1) + 3 N_{\mathrm{def}}(F)\log_2 n$.
\end{proof}

 (L2) is where condensation, at $3.86$, enters. (L3) is where reconstruction enters, and
the two coincide for $k = 3$.

\section*{The decimation phase diagram for \texorpdfstring{$k = 3$}{k = 3}}\label{sec:cavity}

This section computes the evidence~(i) of Section~``\nameref{sec:conditional}''. It is a numerical
cavity computation, not a proof.

Draw $F \sim F_n(\alpha)$ and a uniform solution $\sigma$, and
reveal each variable's value $\sigma_v$ independently with probability $\theta$. The residual
measure is uniform over the solutions that agree with the revealed values. As $\theta$ goes from
$0$ to $1$ this traces the \emph{ideal decimation trajectory} of~\cite{RTS09}. The trial of
Section~``\nameref{sec:conditional}'' follows the same trajectory to the extent that belief propagation
gives the true marginals, so the trajectory is the right object to analyse. Following~\cite{RTS09},
we ask at which $(\alpha, \theta)$ the residual measure has \emph{point-to-set reconstruction}:
whether the values far from a variable, together with the revealed ones, carry information about
it that the revealed ones alone do not.

The cavity method predicts that belief propagation gives asymptotically exact marginals as long as
the measure is not condensed~\cite{KMRSZ07}. Write $\alpha^*$ for the density below which the
trajectory never meets condensation, and $\alpha^*_{\mathrm{rec}}$ for the density below which it
never meets reconstruction. Reconstruction appears no later than condensation, so
$\alpha^*_{\mathrm{rec}} \le \alpha^*$. Above the condensation threshold $\alpha_c$ the measure is
already condensed at $\theta = 0$, so $\alpha^* \le \alpha_c$. Absence of reconstruction is also
the condition under which belief propagation started from uninformative messages is expected to
reach the right fixed point, which the trial needs.

For $k = 3$ the dynamical and condensation thresholds of the undecimated measure coincide,
$\alpha_d = \alpha_c \approx 3.86$~\cite{KMRSZ07}, so there is no clustered, non-condensed phase.
The work~\cite{RTS09} computed the decimation diagram for $k = 4$, where reconstruction appears at
$\theta > 0$ from $\alpha \approx 9.05$, below $\alpha_d \approx 9.38$, and for $k = 5$. It treated
$k = 3$ as a special case, without a regime of positive complexity, and~\cite{MRS07} observed that
the ideal decimation curve of $k = 3$ has no singularity below $3.86$. What we add is a direct
check of the point-to-set reconstruction gap along the trajectory, with a positive control.

We use population dynamics on the cavity tree of random $k$-SAT, in
which a variable has $\mathrm{Poisson}(k\alpha)$ clauses besides its parent clause, each with
$k - 1$ new variables and independent uniform signs. A population element stands for one cavity
tree, rooted at a variable, and stores
\begin{itemize}
  \item $p$, the replica-symmetric cavity marginal $P(x = 1)$ of the undecimated measure;
  \item for each root value $y \in \{0,1\}$, a number $f[y]$: the root's posterior
        $P(x = 1)$ given the revealed variables of a tree drawn with root value $y$;
  \item for each $y$, a number $b[y]$: the same posterior given, in addition, the values of
        all variables at the tree's far boundary.
\end{itemize}
The population of $p$ is first brought to its fixed point. Then $f[y] \gets p$ and
$b[y] \gets y$, the uninformative and the fully informative boundary conditions
of~\cite{MM06}. A sweep replaces $N$ uniformly random elements by new ones. A new element draws
its clauses, with children taken uniformly from the population, and updates $p$ by belief
propagation. For each $y$, it then:
\begin{enumerate}[label=(\arabic*)]
  \item draws the children's values from the tree measure given $x = y$: independent
        $\Ber(p_j)$ if $y$ satisfies the clause, and the same product conditioned on the clause
        being satisfied otherwise;
  \item reveals each child with probability $\theta$; a revealed child sends a hard field at its
        value, and an unrevealed child $j$ with value $x_j$ sends $f_j[x_j]$ (respectively
        $b_j[x_j]$);
  \item sets $f[y]$ and $b[y]$ to the belief-propagation posteriors computed from these
        inputs.
\end{enumerate}
After $t$ sweeps, $b$ carries the boundary information of depth roughly $t$. The
\emph{gap}
\[
  g \coloneqq \mathbb{E}_i \sum_{y \in \{0,1\}} p_i(y)\, \bigl|b_i[y] - f_i[y]\bigr|,
  \qquad p_i(1) = p_i,\ p_i(0) = 1 - p_i,
\]
measures what the boundary adds. There is reconstruction at $(\alpha, \theta)$ when $g$ stays
positive as $t \to \infty$. Because $f$ and $b$ are computed from the same sampled children,
values and reveal coins, their difference contracts to machine precision when there is no
reconstruction: the finite population sets no noise floor on $g$.

\begin{figure}[!htb]
\centering
\includegraphics{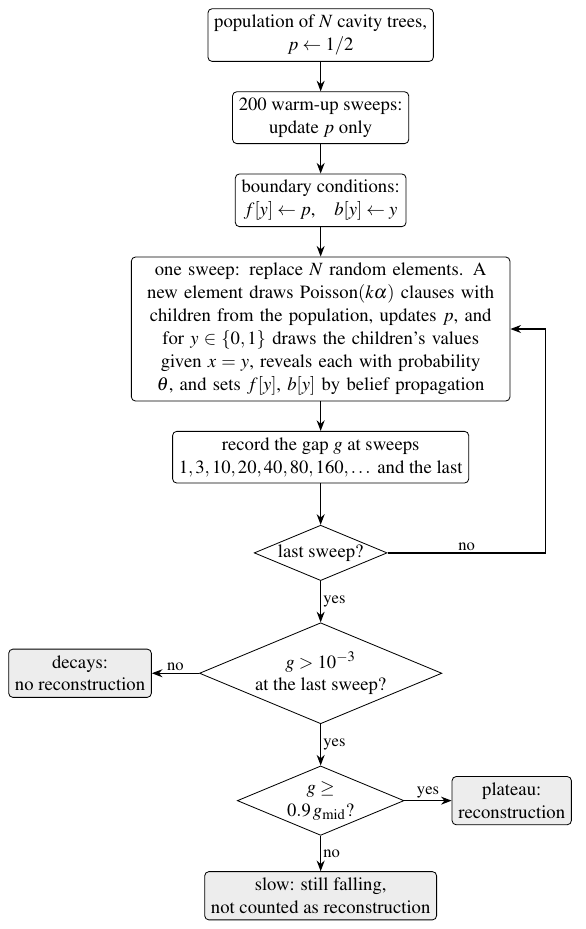}
\caption{\textbf{The cavity experiment for one point $(\alpha, \theta)$.}
The program is \texttt{popdyn.cpp}, and the three classes are those of Fig~\ref{fig:phase}. Here $g_{\mathrm{mid}}$ is the gap at the last recorded
sweep in the first half of the run: sweep $160$ for runs of $400$ or $600$ sweeps, $320$ for runs
of $800$.}
\label{fig:flow-cavity}
\end{figure}

Fig~\ref{fig:flow-cavity} shows how one point $(\alpha, \theta)$ is computed,
and Table~\ref{tab:cavity-runs} lists the five grids. Each point is one population with seed~$1$
and $200$ warm-up sweeps for $p$; the robustness runs below add seeds $2$ and $3$ and a larger
population.

\begin{table}[!ht]
\centering
\small
\caption{\textbf{The five grids of the cavity experiment.}}
\label{tab:cavity-runs}
\begin{tabular}{lcrrll}
\toprule
grid & $k$ & $N$ & sweeps & $\alpha$ & $\theta$ \\
\midrule
coarse & $3$ & $40\,000$ & $400$ & $3.52, 3.65, 3.75, 3.80,$ & $0, 0.02, 0.05, 0.1, 0.15,$ \\
       &     &          &       & $3.86, 3.90, 4.00, 4.20$ & $0.2, 0.3, 0.4, 0.5, 0.7$ \\
\rowsep
fine   & $3$ & $80\,000$ & $800$ & $3.80, 3.83, 3.85, 3.86$ & $0, 0.0025, 0.005, 0.01, 0.015,$ \\
       &     &          &       &                          & $0.03, 0.04, 0.06, 0.08$ \\
\rowsep
above $\alpha_d$ & $3$ & $40\,000$ & $800$ & $3.90, 4.00, 4.20$ & $0, 0.0005, 0.001, 0.002,$ \\
       &     &          &       &                          & $0.005, 0.01$ \\
\rowsep
control & $4$ & $20\,000$ & $400$ & $8.8, 9.0, 9.1,$ & $0, 0.02, 0.04, 0.06,$ \\
        &     &          &       & $9.2, 9.3, 9.45$ & $0.08, 0.1, 0.15$ \\
\rowsep
control, small $\theta$ & $4$ & $20\,000$ & $400$ & $9.45, 9.6$ & $0, 0.0005, 0.001, 0.002,$ \\
        &     &          &       &                  & $0.005, 0.01$ \\
\bottomrule
\end{tabular}
\begin{flushleft}Every point has $200$ warm-up sweeps and seed
$1$. The robustness runs are listed at the end of this section.
\end{flushleft}
\end{table}

\begin{figure}[!htb]
\centering
\includegraphics{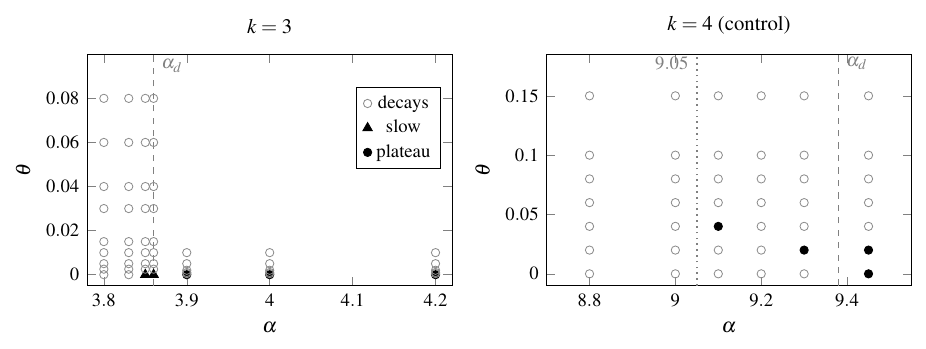}
\caption{\textbf{Where the gap $g$ survives.}
Each point is one population; filled discs are plateaus,
triangles decay slowly ($g > 10^{-3}$ at the last sweep but falling), open circles decay. Left:
$k = 3$, fine grid and the runs above $\alpha_d$ ($800$ sweeps); the coarse grid of the text, up to
$\theta = 0.7$, adds no plateau at $\theta > 0$. Right: $k = 4$ ($400$ sweeps); reconstruction appears at
$\theta > 0$ from $\alpha \approx 9.05$, as in~\cite{RTS09}.}
\label{fig:phase}
\end{figure}

As a positive control we ran the same code at $k = 4$, where it reproduces
the diagram of~\cite{RTS09} (Table~\ref{tab:k4} and Fig~\ref{fig:phase}, right). At $\theta = 0$ the gap decays up to
$\alpha = 9.3$ and has a plateau at $9.45$, on either side of $\alpha_d \approx 9.38$. At
$\theta > 0$ it finds stable plateaus below $\alpha_d$, at $(9.1, 0.04)$ and $(9.3, 0.02)$, and
none at $8.8$ or $9.0$. So the left edge of the reconstruction region at $\theta > 0$ lies
between $9.0$ and $9.1$, consistent with the critical point
$(\alpha^*, \theta^*) \approx (9.05, 0.045)$ of~\cite{RTS09}, near which we find the plateau at
$(9.1, 0.04)$. The region is narrow in $\theta$: at
$\alpha = 9.2$ it falls between the grid points $0.02$ and $0.04$, where $g$ decays to below
$10^{-7}$ by $400$ sweeps. Above $\alpha_d$, at small $\theta$, reconstruction survives: at
$\alpha = 9.45$ and at $\alpha = 9.6$, which is above the condensation threshold $\alpha_c \approx
9.55$ of $k = 4$, the gap keeps a plateau of $0.29$--$0.32$ at every
$\theta \in \{0.0005, 0.001, 0.002, 0.005, 0.01\}$.

\begin{table}[!ht]
\centering
\caption{\textbf{The control at $k = 4$: the gap $g$ after $400$ sweeps.}}
\label{tab:k4}
\begin{tabular}{lcccc}
  \toprule
  $\alpha \backslash \theta$ & $0$ & $0.02$ & $0.04$ & $0.06$ \\
  \midrule
  $8.8$  & $2\cdot10^{-16}$ & $2\cdot10^{-16}$ & $8\cdot10^{-12}$ & $3\cdot10^{-5}$ \\
  \rowsep
  $9.0$  & $2\cdot10^{-16}$ & $2\cdot10^{-13}$ & $3\cdot10^{-6}$  & $3\cdot10^{-7}$ \\
  \rowsep
  $9.1$  & $2\cdot10^{-16}$ & $3\cdot10^{-11}$ & $\boldsymbol{0.17}$  & $1\cdot10^{-10}$ \\
  \rowsep
  $9.2$  & $3\cdot10^{-15}$ & $8\cdot10^{-8}$  & $3\cdot10^{-8}$  & $2\cdot10^{-14}$ \\
  \rowsep
  $9.3$  & $2\cdot10^{-13}$ & $\boldsymbol{0.26}$  & $3\cdot10^{-12}$ & $2\cdot10^{-16}$ \\
  \rowsep
  $9.45$ & $\boldsymbol{0.31}$  & $\boldsymbol{0.24}$  & $5\cdot10^{-16}$ & $3\cdot10^{-17}$ \\
  \bottomrule
\end{tabular}
\begin{flushleft}Bold entries are plateaus (above
$10^{-3}$, and at least $90\%$ of the value at sweep $160$; Fig~\ref{fig:flow-cavity}). At $\theta \ge 0.08$ the gap decays at every
$\alpha$.
\end{flushleft}
\end{table}

\begin{figure}[!htb]
\centering
\includegraphics{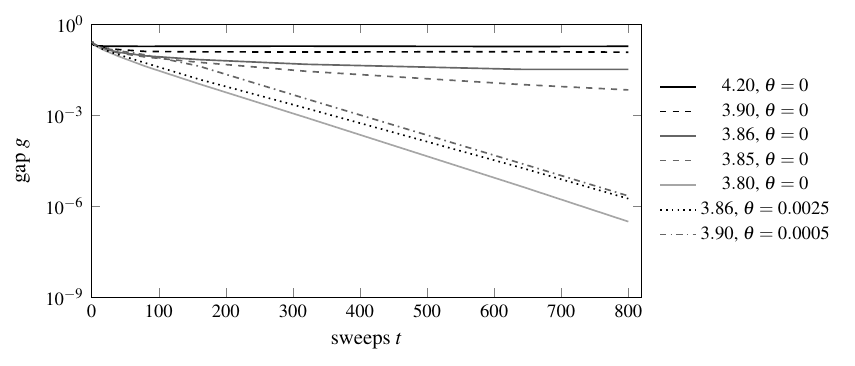}
\caption{\textbf{The gap against the number of sweeps for $k = 3$.}
At $\theta = 0$ it has a plateau above
$\alpha_d$, is at or near the transition at $\alpha_d$, and decays below it. Revealing a small fraction of a
solution makes it decay at every density shown.}
\label{fig:decay}
\end{figure}

At $\theta = 0$ the recursion is the usual reconstruction
recursion of the undecimated measure, which checks the code against the known
$\alpha_d \approx 3.86$ (Table~\ref{tab:theta0}). The gap decays for $\alpha \le 3.80$ and has a
plateau for $\alpha \ge 3.90$. At $3.86$ it is at or near the transition: it falls from $0.09$ at
$80$ sweeps to $0.033$ at $640$, and is nearly flat from there to $800$ (Fig~\ref{fig:decay}),
so it is classed as slow. The decay is slow even well inside the
window: $g$ shrinks by a factor $e$ about every $12$ sweeps at $\alpha = 3.52$ and every $60$
sweeps at $\alpha = 3.80$.

\begin{table}[!ht]
\centering
\caption{\textbf{The gap $g$ after $t$ sweeps at $\theta = 0$, $k = 3$ (coarse grid).}}
\label{tab:theta0}
\begin{tabular}{lcccccc}
  \toprule
  $\alpha$ & $3.52$ & $3.75$ & $3.80$ & $3.86$ & $3.90$ & $4.20$ \\
  \midrule
  $t = 80$  & $4\cdot10^{-4}$  & $2\cdot10^{-2}$ & $4\cdot10^{-2}$ & $9\cdot10^{-2}$ & $0.127$ & $0.187$ \\
  \rowsep
  $t = 160$ & $5\cdot10^{-7}$  & $1\cdot10^{-3}$ & $1\cdot10^{-2}$ & $7\cdot10^{-2}$ & $0.123$ & $0.189$ \\
  \rowsep
  $t = 400$ & $7\cdot10^{-16}$ & $1\cdot10^{-6}$ & $2\cdot10^{-4}$ & $4\cdot10^{-2}$ & $0.124$ & $0.187$ \\
  \bottomrule
\end{tabular}
\end{table}

For $\theta > 0$ we found no reconstruction at any point of any $k = 3$ grid.
\begin{itemize}
  \item On the coarse grid, $g$ decays at every $\theta \ge 0.02$ and every $\alpha$, and is $0$ to
        machine precision for $\theta \ge 0.4$. At $\theta = 0.02$ and $t = 400$ it is
        $5\cdot10^{-10}$ at $\alpha = 3.86$ and $2\cdot10^{-15}$ at $\alpha = 4.0$. The decay
        speeds up with $\theta$: at $\alpha = 3.80$ and $t = 80$, $g$ is $4\cdot10^{-2}$,
        $5\cdot10^{-3}$, $2\cdot10^{-4}$ and $5\cdot10^{-6}$ for $\theta = 0$, $0.02$, $0.05$
        and $0.1$.
  \item On the fine grid, $g$ decays at every $\theta > 0$ and every $\alpha$, to at most
        $1.8\cdot10^{-6}$ after $800$ sweeps (the largest values are at $\theta = 0.0025$) and to
        machine precision for $\theta \ge 0.015$. Table~\ref{tab:fine} shows the two values of
        $\alpha$ closest to $\alpha_d$.
\end{itemize}

\begin{table}[!ht]
\centering
\caption{\textbf{The gap $g$ after $800$ sweeps near $\theta = 0$, $k = 3$ (fine grid).}}
\label{tab:fine}
\begin{tabular}{lcccccc}
  \toprule
  $\theta$ & $0$ & $0.0025$ & $0.005$ & $0.01$ & $0.015$ & $0.03$ \\
  \midrule
  $\alpha = 3.85$ & $7\cdot10^{-3}$ & $2\cdot10^{-6}$ & $1\cdot10^{-8}$ & $3\cdot10^{-12}$ & $1\cdot10^{-15}$ & $3\cdot10^{-16}$ \\
  \rowsep
  $\alpha = 3.86$ & $3\cdot10^{-2}$ & $2\cdot10^{-6}$ & $6\cdot10^{-9}$ & $9\cdot10^{-13}$ & $9\cdot10^{-16}$ & $3\cdot10^{-16}$ \\
  \bottomrule
\end{tabular}
\end{table}

So, on these grids, the ideal decimation trajectory of random 3-SAT stays out of reconstruction
for every $\alpha < \alpha_d \approx 3.86$: $\alpha^*_{\mathrm{rec}} \approx 3.86$. Since
$\alpha^*_{\mathrm{rec}} \le \alpha^* \le \alpha_c \approx 3.86$, this also places $\alpha^*$ at
about $3.86$. Unlike $k = 4$, there is no gap between $\alpha^*$ and $\alpha_d$: revealing even a
$0.25\%$ fraction of a uniform solution removes the slow decay seen at $\theta = 0$.

Above $\alpha_d$, the gap also decays at every $\theta > 0$ we tried, down to $\theta = 0.0005$
at $\alpha \in \{3.90, 4.00, 4.20\}$ ($g \le 5\cdot10^{-6}$ after $800$ sweeps; three seeds at
$3.90$ and $4.20$, one at $4.00$),
while $\theta = 0$ has a plateau. We do not have an explanation. We first guessed that, with no
clustered, non-condensed phase, revealing a small fraction of a solution pins the dominant cluster
and leaves a measure that does not reconstruct. The $k = 4$ control contradicts that guess: above
its condensation threshold, reconstruction survives the same small fractions. A possible lead is
that at $k = 3$ the transition at $3.86$ is a local instability~\cite{KMRSZ07,BMW00}, unlike the
transition at $k = 4$, and a small revealed fraction may be enough to remove a local instability;
we have not tested it. These rows say
nothing about the trial at $\theta$ close to $0$ either way: above $\alpha_c$ the tree measure is
no longer the law of a uniform solution.

We repeated the key points with seeds $2$ and $3$ ($N = 80\,000$, $800$
sweeps) and with $N = 160\,000$ ($600$ sweeps).
\begin{itemize}
  \item \emph{Seeds.} The points are $\alpha \in \{3.80, 3.85, 3.86\}$ with
        $\theta \in \{0, 0.0025, 0.01\}$, and $\alpha \in \{3.90, 4.20\}$ with
        $\theta \in \{0, 0.0005, 0.002\}$. Every point gets the same class with seeds $2$ and $3$. The plateaus agree
        within $5\%$ ($0.124$ and $0.122$ against $0.120$ at $3.90$; $0.186$ and $0.187$ against
        $0.188$ at $4.20$). The slow points agree ($0.034$ and $0.034$ against $0.03$ at $3.86$),
        and the decaying points agree within a factor $2$.
  \item \emph{Population size.} The points are $(3.85, 0.0025)$, $(3.86, 0)$, $(3.90, 0.0005)$
        and $(4.20, 0)$. At $N = 160\,000$ every point gets the same class, and the decay rates
        over the second half of the sweeps differ by $17$--$22\%$, over different sweep windows.
\end{itemize}

The remaining limits of the computation are finite $N$, finitely many sweeps, and finite grids. The control shows that a reconstruction region at $\theta > 0$ can be
narrow, about $0.02$ wide in $\theta$ near its tip, so a narrower region at $k = 3$ could fall
between grid points; the fine grid near $\theta = 0$ is spaced by $0.0025$ to $0.02$, the coarse
grid by $0.02$ to $0.2$. The slow decay at $\theta = 0$ means the trial's first steps see
correlations over many levels of the tree, more than the radius $\sim \log n / \log(6\alpha)$ up to
which a random formula looks like a tree.

\section*{Conclusion and future work}\label{sec:future}

\emph{Theoretical findings.} A Las Vegas algorithm that sets variables one at a time from estimated
conditional probabilities succeeds with probability at least $2^{-K}$, where $K$ is the
information-theoretic cost of the estimates' errors along the run of a uniform solution. The bound
holds for any sound closure and any estimator, randomised ones included, and over finite domains.
Every output comes with the exact probability with which it was produced, so the importance weights
give an unbiased estimate of the number of solutions when every solution can be produced, and a
certified lower bound on it in every case. The same probabilities bound the distance of the outputs
from uniform from below, without any exact count. That distance is the divergence in the opposite
direction to $K$: a small $K$ is sufficient for a high success rate but not necessary, and a sampler
can trade uniformity for success. A $K$ that is logarithmic in the number of variables holds on every
satisfiable formula if and only if $\mathrm{NP} = \mathrm{RP}$. For random 3-SAT, an estimator whose $K$ grows logarithmically in
probability would give a Las Vegas algorithm that runs in polynomial time on all but any fixed
fraction of formulas (Corollary~\ref{cor:main}), beyond the densities at which any algorithm is
proved to succeed with high probability. We conjecture that one exists at every density below the
satisfiability threshold (Conjecture~\ref{conj:main}). Below $3.86$ our candidate is belief
propagation, supported by the cavity computation of Section~``\nameref{sec:cavity}'', and
Section~``\nameref{sec:towards}'' reduces the hypothesis to three lemmas with no known proof. Above $3.86$ we
have no evidence either way, and state the conjecture there because nothing known rules it out.

\emph{Practical findings.} On timetabling, the choice of estimator decides how well the sampler
works. On generated instances, belief propagation on the encoding is confidently wrong around the
deterministic loops formed by the auxiliary variables, while belief propagation on the timetabling
problem itself has a cost $K$ of at most about one bit, estimates the number of timetables within
$2\%$ and the flexibility map within $0.3$ percentage points, and samples close to uniformly. On a
real competition instance with $676$ classes and about three million forbidden pairs, rerunning
only the failed runs delivered every one of $72$ requested timetables with each of two estimators,
and every output passed a validator. The work per step is set by the structure of the
problem, never by a time limit or a test of convergence. Belief propagation succeeds in $31\%$ of
runs; belief propagation steered by a
Markov chain's witness succeeds in $65\%$, with outputs at least $327$ bits from uniform. The first
is the one to use when the outputs must represent the feasible timetables as a whole, the second
when the aim is many different feasible timetables quickly. The instance has more than $2^{2165}$
feasible timetables, with probability at least $0.99$. On it, the importance weights of both
estimators are too uneven to estimate the number of timetables or the flexibility map, which is
the main limitation of the application as it stands.

\emph{Further work on the application.} A comparison with the samplers cited in
Section~``\nameref{sec:intro}'', CMSGen, UniGen and QuickSampler, on speed, diversity and uniformity, with
Barbarik as the test of uniformity. More competition instances, including those whose required
constraints, such as a maximum daily load, are not between pairs of classes, and need factors over
more than two classes. And an estimator that uses the witness's certificate that the current state is not a dead end, without the witness
dominating the draw as it does in the product of Section~``\nameref{sec:itc}''. 

\emph{Further work on NP versus RP.} By Proposition~\ref{prop:worst} and
Remark~\ref{rem:refute}, the conjecture is tied to $\mathrm{NP}$ versus $\mathrm{RP}$ in both
directions, and no partial result on either side implies either. Information alone does not settle
it, since the exact conditional probabilities give $K = 0$ on every satisfiable formula. The question
is whether they can be computed in polynomial time. What would bring it closer is a property of the
formulas that bears on this, together with a way to check where the property holds. Our candidate is
point-to-set reconstruction along the decimation trajectory, under which the conditional probability
of a variable depends on values far from it. The overlap gap property plays a similar role for stable
algorithms~\cite{GS17,BH21}. The further work is to determine whether reconstruction along the
trajectory prevents approximate counting of the restricted formulas, and where it occurs. Section~``\nameref{sec:cavity}'' gives numerical evidence on the second
part: on our grids, reconstruction does not occur below $3.86$. Above $3.86$ the computation does
not settle the question.

\end{document}